\documentclass[11pt,notitlepage,twoside]{article}
\usepackage{latexsym}
\usepackage{amsmath}
\usepackage{amsfonts}
\usepackage{amssymb}
\usepackage{enumerate}
\usepackage[latin1]{inputenc} 
\renewcommand{\a }{\alpha }
\renewcommand{\b }{\beta }
\renewcommand{\d}{\delta }

\newcommand{\D }{\Delta }

\newcommand{\e }{\varepsilon }
\newcommand{\eps }{\varepsilon }
\newcommand{\g }{\gamma}

\renewcommand{\l }{\lambda }
\renewcommand{\L }{\Lambda }

\newcommand{\n }{\nabla }

\renewcommand{\O }{\Omega }

\newcommand{\ov}{\overline}

\newcommand{\be}{\begin{equation}}
\newcommand{\ee}{\end{equation}}
\newenvironment{pf}{\noindent{\bf Proof.}\enspace}{
\hfill$\Box$\medskip}

\newcommand{\R}{\mathbb{R}}
 
\newtheorem{thm}{Theorem}[section]
\newtheorem{pro}[thm]{Proposition}
\newtheorem{lem}[thm]{Lemma}

\numberwithin{equation}{section} 

\title { \Large \textbf{Multispike Solutions for a slightly subcritical \\
elliptic problem with non-power nonlinearity }}
\author{{\bf\large Mohamed Ben Ayed}
{\bf\large}\vspace{1mm}\\
 {\it\small Department of Mathematics,  College of Science, Qassim University, Saudi Arabia,}\\
{\it\small University of Sfax, Faculty of Sciences of Sfax,
Tunisia}\\
{\it\small e-mail: Mohamed.Benayed@fss.rnu.tn}\vspace{1mm}\\
{\bf\large Habib Fourti} 
{\bf\large}\vspace{1mm}\\
{\it\small University of Monastir, Faculty of Sciences of
Monastir, Tunisia},\\ {\it\small Laboratory LR 13 ES 21, University of Sfax, Faculty of Sciences of Sfax, Tunisia}\\
{\it\small e-mail: habib40@hotmail.fr, habib.fourti@fsm.rnu.tn}
\\
\vspace{1mm}{\bf\large Rabeh Ghoudi}
{\bf\large}\vspace{1mm}\\
{\it\small University of Gabes, Faculty of Sciences of Gabes,
Tunisia}\\
{\it\small e-mail: ghoudi.rabeh@yahoo.fr}\vspace{1mm}}

\begin{document}

\date{}

\maketitle

{\footnotesize \noindent {\bf Abstract.} In this paper, we are
concerned with the following elliptic equation
$$\left\{\begin{array}{rrl}-\Delta u&=& |u|^{4/(n-2)}u/[\ln
(e+|u|)]^\e\hbox{ in } \Omega,\\ u&=&0 \hbox{ on }\partial \Omega,
\end{array}
\right.$$ where $\Omega $
 is a smooth bounded open domain in $\mathbb{R}^n, \ n\geq 3$ and
 $\e>0$.
Clapp et al. in Journal of Diff. Eq. (Vol 275) proved that there
exists a
single-peak positive solution for small $\e$ if $n \geq 4$. \\
Here we construct positive as well as changing sign solutions
concentrated at several points at the same time.}

{\footnotesize\section{Introduction and results}
 Let us consider the
nonlinear elliptic problem:
\begin{eqnarray*} (P_{\e})\qquad \qquad  \begin{cases}
 -\Delta u= f_{\e}(u)    & \mbox{in} \, \, \O,  \\
 \quad \ \ u= 0     &  \mbox{on}\, \, \partial \O,
\end{cases}
\end{eqnarray*}
where $\O$ is a smooth bounded domain in $\R^n$, $n\geq 3$ and
\begin{equation}\label{nonpower}
f_\e(u):=\frac{|u|^{p-1}u}{[\ln (e+|u|)]^\e}, \ \e\geq 0, \
p=\frac{n+2}{n-2}.
\end{equation}
Here  $p+1=2n /(n-2)$ is the critical Sobolev exponent for the
embedding of $ H^1_0({\O})$ into $L^{p+1}(\O)$. \\Recently, the
nonlinearity $f_\e$ has attracted a lot of attention (see for
instance \cite{Pistoianew},\cite{Har},\cite{Pardo}). When $\e=0$, we
get $f_0(u)=|u|^{p-1}u$. In this case the variational problem
corresponding to $(P_0)$ happens to be lacking of compactness and
this is due to the noncompactness of the embedding of $ H^1_0({\O})$
into
$L^{p+1}(\O)$.\\
Despite that  the nonlinearity $f_\e$ is very close to the critical
growth, problem $(P_\e)$ is considered as a slightly subcritical
one. Indeed, the nonlinearity $f_\e $ satisfies the condition
\begin{equation}\label{harcondition}\lim_{s\rightarrow \infty} \frac{f (x, s) }{|s|^
\frac{n+2}{n-2}} = 0\hbox{ uniformly with respect to } x \in \Omega
,\end{equation} since $f_\e$ is defined by \eqref{nonpower} and
independent of the variable $x$. This new assumption, recently
introduced by Harrabi in \cite{Har} guarantees the existence of
solutions for problem $(P_\e).$ In fact, the Euler Lagrange
functional associated to $(P_\e)$ satisfies Palais Smale condition
and compactness is
recovered thanks to \eqref{harcondition}. 
This is in some sense very similar to what happens in the
subcritical regime where we usually use Ambrosetti-Rabinowitz
conditions and the following subcritical polynomial growth condition
$$f(x,s)\leq C (|s|^p+1), \hbox{ for all }(x,s)\in \Omega\times \mathbb{R}\
\hbox{ for some }p\in [1, \frac{n+2}{n-2})$$
to ensure the compactness of a bounded Palais Smale sequence. For more details of these aspects we refer to \cite{Har}.
Recently, Clapp et al. constructed in \cite{Pistoianew} a
single-peak solution for problem $(P_\e)$ and they proved that any
$x_0$ non-degenerate critical point of the Robin function generates
a family of solutions concentrating around $x_0$ as $\e$ goes to
$0$. In \cite{Pardo}, the authors have analyzed the asymptotic
behavior of radially symmetric solutions of $(P_\e)$ when $\Omega$
is a ball. The result of Clapp et al. \cite{Pistoianew} was a first
step towards establishing the existence of blowing up solutions to
problem $(P_\e)$. In this direction, our
 main result provides the existence of positive as well as changing sign solutions that blow up and/or blow down at
different points in $\Omega$. In fact, problem $(P_\e)$ shares many
aspects with the following nonlinear subcritical elliptic problem
\begin{eqnarray}\label{subpb}\begin{cases}
 -\Delta u= |u|^{p-1-\e}u    & \mbox{in} \, \, \O,  \\
 \quad \ \ u= 0     &  \mbox{on}\, \, \partial \O,
\end{cases}
\end{eqnarray}
where $\e>0$. Thus one can expect the extension of existence type
results obtained for \eqref{subpb} to our problem $(P_\e)$ as
conjectured in \cite{Pistoianew}. But  the non-power nonlinearity
\eqref{nonpower} gives rise to some technical difficulties. In this
work, we were able to overcome these issues by some careful estimates.\\
Problem \eqref{subpb} has been investigated widely in the last
decades, see for example \cite{BLR,BMP,Han,MP2,MP1,PW,Pistoia,R,R2}.
In the sequel, we review some known facts related to \eqref{subpb}.
In \cite{Han} and \cite{R2}, the asymptotic behavior of least energy
solutions (hence positive) was analyzed as $\e\rightarrow 0$. It was
proved that, a single peak solution blows up and concentrates at a
critical point of the Robin's function as $\e$ goes to $0$.  This
study was completed in \cite{BLR} in the case of general positive
blowing up solutions. Conversely, by a finite reduction procedure,
Bahri, Li and Rey proved the existence of a family of positive
solutions of $\eqref{subpb}$ concentrating at several points of
$\Omega.$ Later, Rey \cite{R3} extended the same results to
dimension $3$. Concerning changing sign solutions, Bartsch,
Micheletti, and Pistoia \cite{BMP} proved that the same method as
Bahri-Li-Rey's produces multispike solutions
 blowing up or
blowing down at different points characterized as solutions of a
system of equations defined explicitly in terms of the gradients of
the Green function and its regular part.\\ 
In \cite{PW}, the authors proved that if $\Omega$ is symmetric with
respect to the $x_1,\ldots, x_n$ axes, problem \eqref{subpb} has a
sign-changing solution with the shape of a tower of bubbles with
alternate signs, centered at the center of symmetry of the domain.
Few years later, Musso and Pistoia \cite{MP1} were able to extend
such a result to a general domain.
When the Laplacian operator in \eqref{subpb} is replaced by the
biharmonic one, the
 existence of multispike and bubble towers changing sign
solutions for the counterpart of \eqref{subpb} have been studied in
\cite{BG} and \cite{DG} respectively.

Going back to $(P_\e)$, the existence of multispike solutions has
not been studied yet and the main purpose of this paper is to focus
on this issue. More precisely, we construct families of solutions
for the equation, which blow-up positively or negatively and
concentrate in $m$
different points of $\Omega$ as the parameter $\e$ goes to $0$.\\
In order to state our main result, we introduce some notations.\\
The space $H^1_0(\O)$ is equipped with the norm $\|.\|$ and its
corresponding inner product $\left<.,.\right>$ defined by
$$ \|u\|^2=\int_{\O}|\n
u|^2;\quad \left<u,v\right>=\int_{\O}\n u \n v,\quad u,v \in
H^1_0(\O).$$ For $a\in \O$ and $\l>0$, let
 \begin{equation}\label{d}
\d_{(a,\l)}(y)=\frac{c_0\l^{(n-2)/2}}{\left(1+\l^2|y-a|^2\right)^{(n-2)/2}},
\hbox{ where } c_0:=(n(n-2))^{(n-2)/4}.
 \end{equation}
The constant $c_0$ is chosen such that $\d_{(a,\l)}$ is the family
of solutions of the following problem
\begin{equation}\label{e}
 -\Delta u=u^{(n+2)/(n-2)},\;\;u>0\; \;\mbox{in}\; \;\R^n.
\end{equation}
Notice that the family $ \d_{(a,\l)}$ achieves the best Sobolev
constant \begin{equation}\label{S_n}S_n:=\inf\{\|\n
u\|^2_{L^2(\R^n)}\|u\|^{-2}_{L^{2n/(n-2)}(\R^n)}:u\not\equiv 0, \n u
\in (L^{2}(\R^n))^n \,\mbox{and}\;u \in L^{2n/(n-2)}(\R^n)\}.
\end{equation}
 We denote
by $P\d_{(a,\l)}$ the projection of $\d_{(a,\l)}$ onto $ H^1_0(\O)$,
defined by
\begin{eqnarray}\label{q1}
-\Delta P\d_{(a,\l)}=-\D \d_{(a,\l)} \mbox{ in } \O, \quad
P\d_{(a,\l)}=0 \mbox{ on}\;
\partial \O.
\end{eqnarray}
We will denote by $G$ the Green's function and by $H$ its regular
part, that is
$$G(x,y)=|x-y|^{2-n}-H(x,y) \quad\mbox{ for }\,\,
(x,y)\in \O^2,$$ and for $x\in\O$, $H$ satisfies
\begin{eqnarray*}   \begin{cases}
\Delta  H(x,.)=0   & \mbox{in} \, \, \O,
\\H(x,y)=|x-y|^{2-n},&\mbox{for}\, \,y\in
\partial \O.
\end{cases}
\end{eqnarray*}
Next we describe the solutions that we are looking for. Let $m$ be
an integer and $(\gamma_1,\ldots,\gamma_m)\in\{-1,1\}^m$, we
construct solutions of the form
\begin{equation}\label{u_eps}
u_\e= \sum_{i=1}^m\alpha_i\gamma_iP\delta_{(a_i,\lambda_i)}+v,
\end{equation}
where $(\a_1,...,\a_m)\in (0,+\infty)^m$, $(\l_1,...,\l_m)\in
(0,+\infty)^m$ and $(a_1,...,a_m)\in \O^{m}$. The term $v$ has to be
thought as a remainder term of lower order. Let
\begin{align}\label{vorthogonal}&E_{(a,\l)}:=
\Big\{v\in H^1_0(\O):\bigg<v,P\d_{i}\bigg>=\left<v,\frac{\partial
P\d_{i}}{\partial \l_i}\right>= \left<v,\frac{\partial
P\d_{i}}{\partial (a_i)_j}\right>=0\;\forall \; 1\leq j\leq
n,\,\forall \,1\leq i\leq m\Big\},
\end{align}
where $P\d_i=P\d_{(a_i,\l_i)}$ and $( a_i )_j$ is the $j^{th}$ component of $ a_i $.\\

For $x=(a_1,\ldots, a_m)\in\O^m$ and
$(\gamma_1,\ldots,\gamma_m)\in\{-1,1\}^m$, we denote by
$M(x)=(m_{ij})_{1\leq i,j\leq m}$ the matrix defined by
\begin{equation}\label{coef}
m_{ii}=H(a_i,a_i) \ ; \ m_{ij}=-\gamma_i\gamma_jG(a_i,a_j)\ , \
i\neq j,
\end{equation}
and by $\rho(x)$ its least eigenvalue. We also define
$$\begin{array}{rll}
{\bf{F}}_{x}:(0,+\infty)^{m}&\longrightarrow&\mathbb{R}\\
\Lambda=(\L_1,\ldots,\L_{m})&\longmapsto& \frac{1}{2}\L M(x)^t\L-\ln
\L_1\ldots\L_{m}.
\end{array}$$
If $\rho(x) > 0$, ${\bf{F}}_x$ is strictly convex on
$(0,+\infty)^{m}$, infinite on the boundary ; so ${\bf{F}}_x$ has in
$(0,+\infty)^{m}$ a unique critical point $\Lambda(x)$, which is a
minimum. On the subset of $\O^m$
$$\rho^+=\{x\in\O^m/ \rho(x)>0\},$$
we define the function
\begin{equation}\label{18} \widetilde{{\bf{F}}}(x) = {\bf{F}}_x(\Lambda(x))=
\frac{m}{2}- \ln \L_1(x)\ldots \L_m(x) ,\end{equation} whose
differential is given by
\begin{equation*}
\widetilde{{\bf{F}}}'(x) = \frac{1}{2}\L(x)M'(x)^t \L(x) = -
\displaystyle \sum_{i=1}^m \frac{\L_i'(x)}{\L_i(x)}.
\end{equation*}
Now, we are able to state the following result:
\begin{thm}\label{t:11} Let $ n\geq 4$ and $(\gamma_1,\ldots, \gamma_m)\in \{-1,1\}^m$. Let $\overline{x}=(\overline{x}_1,\ldots,\overline{x}_m)\in \rho^+$
be a non-degenerate critical point of $\widetilde{{\bf{F}}}$. Then,
there exists $\e_0>0$ such that for each $\e\in (0,\e_0)$, problem
$(P_\e)$ has a solution $(u_\e)$  of the form:
\begin{equation} \label{h:180}
u_\e= \sum_{i=1}^m\alpha_{i,\e}\gamma_iP\delta_{( a_{i,\e}
,\lambda_{i,\e})}+v_\e,
\end{equation}
where, as $\e\rightarrow 0$,
$$\|v\|\rightarrow 0,\ \left|\a_{i,\e}-1\right|=O(\e\ln|\ln\e|),
\ \lambda_{i,\e}^{-\frac{n-2}{2}}\left(\frac{|\ln
\e|}{\e}\right)^\frac{1}{2}\rightarrow\overline{c}\Lambda_i(\overline{x}),\
|a_{i,\e}-\overline{x}_i|=O\left(\frac{\ln|\ln\e|}{|\ln\e|}\right).$$
Here, $\overline{c}$ is some positive constant. Moreover, we have
$|\nabla u_{\eps}|^{2}\rightharpoonup S_n^{n/2}\displaystyle
\sum^{\ell}_{i=1}\delta_{\overline{x}_i}$ in $\mathcal{D}'(\Omega)$
when $\eps\rightarrow 0$, where $\d_y$ denotes the Dirac mass at the
point $y$.
\end{thm}
Let us remark that the solutions of $(P_\e)$ we obtained in this
paper do not converge to solutions of the limit problem $(P_0)$
(even if $(P_0)$ has solutions) since
they converge weakly to zero as $\e\rightarrow 0$.\\
We think that the same result holds true even for $n=3$. The
restriction on the dimension in our result comes from the bad
estimate of the $v$-part of the solution in dimension three. To
handle the case of three dimensional bounded domains, one can
proceed as in \cite{R3} and that by improving the odd part of $v$.\\
Our result Theorem \ref{t:11} extends both existence results in
\cite{BLR} and
\cite{BMP} to a non-power nonlinearity.\\
Taking $\gamma_i=1$ for each $1\leq i\leq m$, a non-degenerate
critical point $\overline{x}=(\overline{x}_1,\ldots,\overline{x}_m)$
of $\widetilde{{\bf{F}}}$ in $\rho^+ $ generates a family of
positive solutions of $(P_\e)$ which blow up and concentrate at
m-different points. Indeed, since $ |u_{\e}^-|_{L^{p+1}}$ is small,
arguing as in \cite{BLR}, the solution $u_{\e}$ has to be positive
in $\O$. Such an assumption was firstly introduced in \cite{BLR} to
study the
subcritical problem \eqref{subpb}. \\
As in \cite{Pistoianew}, the blow up rate of our solutions $u_\e$
satisfies $\|u_\e\|_\infty\sim c (|\ln \e|\e^{-1})^{1/2}$ for both
positive and changing sign solutions as $\e$ goes to $0$. Compared
to multispike solutions of \eqref{subpb}, solutions given by Theorem
\ref{t:11} blow up faster. In fact, as described above, the
concentration speed $\l_i$'s are of order $(|\ln \e
|\e^{-1})^{1/(n-2)}$ which differ from the usual power nonlinearity
as stated in \cite{BLR} and \cite{BMP} where the $\l_i$'s are of
order $\e^{-1/(n-2)}$. Our choice will be justified by the expansion
of $\left< \nabla I_\e(u),\l_i\partial P\d_i/\partial \l_i\right>$
given in Proposition \ref{p23} where $I_\e$ is
 the functional associated to $(P_\e)$ introduced in \eqref{g0}.\\
The proof of our result is based on the finite reduction method
introduced in \cite{BLR}. This method has been widely used to study
elliptic problems involving critical Sobolev exponent with small
perturbations. The proofs of all existence results of blowing up or
blowing down solutions  as the parameter $\e$ goes to zero,
mentioned
before,  rely on this technique. 
We describe in the following our proof.
Firstly, we show that given $\underline{u}=\sum^m_{i=1}
\alpha_i\gamma_i P\delta_i$ , the functional $I_\e (\underline{u} +
v)$ (introduced in \eqref{g0}) has a unique minimum with respect to
the variable $v \in E_{a,\l}$.
\\
This leads to a finite dimensional reduced problem depending of the
variables $\a_i$, the concentration rates $\l_i$ and the
concentration points $a_i$. Our analysis requires a careful
expansion of the gradient of the energy functional $I_\e$. This will
be developed in Section 2. We mention that we follow some ideas in
\cite{Pistoianew} and we also improve some estimates and complete
the analysis started in this work by taking into account the
different interactions among the bubbles which depend of their
respective signs.\\
Then by using a suitable change of variables, we obtain a system
satisfied by the new variables. This system is explicitly given in
terms of the geometry of the domain, namely Green and Robin
functions applied to $\overline{x}$.\\
Lastly, the fact that $\overline{x}$ is a non-degenerate critical
point of
$\widetilde{{\bf{F}}}$ allows us to conclude with a fixed point theorem.\\
Finally, we mention that our analysis of the gradient which is
developed in Section 2 was done in a general setting.  We think that
 the obtained expansions will be useful to describe the blow up profile of positive and changing sign
 solutions, as $\e$ goes to zero.\\
We also think that combining these
 expansions with some ideas in \cite{MP1}, one can construct bubble
 tower solutions. In fact, after the choice of
suitable concentration rates and concentration points, the author in
\cite{MP1} applied similar finite dimensional reduction methods. We
expect that similar arguments produce these kind of solutions for
the problem $(P_\e)$.
\\
\\
The paper is organized as follows: Section 2 is devoted to the
technical framework which includes some asymptotic expansions of the
gradient of the energy functional. In Section 3, we prove our main
result.\\
Throughout this paper, we use the same $c$ to denote various generic
positive constants independent of $\e$.

\section{The Technical Framework}
We introduce the general setting. For $ \e> 0$, we define on $
H^1_0(\O)$ the functional
\begin{align}\label{g0}&I_\e(u)=\frac{1}{2}\int_{\O}|\n
u|^2-\int_{\O}F_\e(u), \hbox{ where }F_\e(s)=\int_0^sf_\e(t)dt.
\end{align}
Note that each critical point of $I_\e$ is a solution of $(P_\e)$.
Since our aim is to construct solutions of the form \eqref{u_eps},
we give the expansion of the gradient of $I_\e$ for $u=
\sum_{i=1}^m\alpha_i\gamma_iP\delta_{(a_i,\lambda_i)}+v$ i.e.  $u$
belongs to a neighborhood of potential concentration sets.
To do so we need some preliminary results. We start with
\cite[Proposition 1]{R}.
\begin{pro}\label{p21} Let $a\in
\O$ and $\l>0$ such that $\l d:=\l d(a, \partial\O)$ is large
enough. For $\varphi_{(a,\l)}=\d_{(a,\l)}-P\d_{(a,\l)}$, we have the
following estimates
\begin{eqnarray*}(a)\;\; 0\leq \varphi_{(a,\l)}\leq \d_{(a,\l)},\quad
(b)\quad \varphi_{(a,\l)}=c_0\frac{H(a,.)}{\l^{(n-2)/2}}+f_{(a,\l)},
\end{eqnarray*}
where $c_0$ is defined in \eqref{d} and $f_{(a,\l)}$ satisfies
\begin{align*}&f_{(a,\l)}=O\biggr(\frac{1}{\l^{\frac{n+2}{2}}d^{n}}\biggr),\quad
\l\frac{\partial f_{(a,\l)}}{\partial
\l}=O\biggr(\frac{1}{\l^{\frac{n+2}{2}}d^{n}}\biggr),
\\&\frac{1}{\l}\frac{\partial f_{(a,\l)}}{\partial
a}=O\biggr(\frac{1}{\l^{\frac{n+4}{2}}d^{n+1}}\biggr).
\end{align*}
\begin{align*}(c)\;&|\varphi_{(a,\l)}|_{2n/(n-2)}=O\biggr(\frac{1}{(\l
d)^{(n-2)/2}}\biggr),\quad
\biggr|\l\frac{\partial\varphi_{(a,\l)}}{\partial
\l}\biggr|_{2n/(n-2)}=O\biggr(\frac{1}{(\l
d)^{(n-2)/2}}\biggr),\\&\|\varphi_{(a,\l)}\|=O\biggr(\frac{1}{(\l
d)^{(n-2)/2}}\biggr),\quad
\biggr|\frac{1}{\l}\frac{\partial\varphi_{(a,\l)}}{\partial
a}\biggr|_{2n/(n-2)}=O\biggr(\frac{1}{(\l d)^{n/2}}\biggr),
\end{align*}
where $|.|_q$ denotes the usual norm in $L^q(\Omega)$ for each
$1\leq q\leq \infty$.
\end{pro}
We also introduce the following lemmas which are used in several
works. For the proof, we refer to \cite{B1}.
\begin{lem}\label{estd}
Let $a\in\Omega$
 and $\lambda > 0$ be such that $\lambda d := \lambda d(a,
 \partial\Omega)$ is very large. There hold
 $$\|P\delta_{(a,\lambda)}\|^2=S_n^{n/2}-\overline{c}_1\frac{H(a,a)}{\l^{n-2}}+O\left(\frac{\ln(\l d)}
 {(\lambda d)^{n}}\right),$$
 $$\left<P\delta_{(a,\lambda)},\lambda\frac{\partial P\delta_{(a,\l)}}{\partial \l}\right>=
 \frac{n-2}{2}\overline{c}_1\frac{H(a,a)}{\l^{n-2}}+O\left(\frac{\ln(\l d)}{(\l d)^n}\right),$$
 $$\left<P\delta_{(a,\lambda)},\frac{1}{\lambda}\frac{\partial P\delta_{(a,\l)}}{\partial a}\right>=
-\overline{c}_1\frac{\frac{\partial H}{\partial
a}(a,a)}{\l^{n-1}}+O\left(\frac{\ln(\l d)}{(\l d)^{n+1}}\right),$$
where $\partial H/\partial a$ denotes the partial derivative of $H$
with respect to the first variable, $S_n$ is introduced in
\eqref{S_n} and
\begin{equation}\label{c_1}\displaystyle
\overline{c}_1:=c_0^{\frac{2n}{n-2}}\int_{\mathbb{R}^n}\frac{1}{(1+|x|^2)^{(n+2)/2}}dx.
\end{equation}
\end{lem}
Note that the constant $S_n$ satisfies
$$S_n^{n/2}=c_0^\frac{2n}{n-2}\int_{\mathbb{R}^n}\frac{1}{(1+|x|^2)^n}
dx.$$
\begin{lem}
Let $a\in\Omega$
 and $\lambda > 0$ be such that $\lambda d := \lambda d(a,
 \partial\Omega)$ is very large. There hold
 $$\displaystyle \int_{\O}P\delta_{(a,\lambda)}^{p+1}=S_n^{n/2}-\frac{2n}{n-2}\overline{c}_1\frac{H(a,a)}{\l^{n-2}}+O\left(\frac{\ln(\l d)}
 {(\lambda d)^{n}}\right)+ \ (\hbox{if }n=3)\ O\left(\frac{1}{(\l d)^2} \right),$$
 $$\displaystyle \int_{\O}P\delta_{(a,\lambda)}^{p}\l\frac{\partial P\delta_{(a,\lambda)}}{\partial \l}=
 2\left<P\delta_{(a,\lambda)},\lambda\frac{\partial P\delta_{(a,\l)}}{\partial \l}\right>+O\left(\frac{\ln(\l d)}
 {(\lambda d)^{n}}\right)+ \ (\hbox{if }n=3)\ O\left(\frac{1}{(\l d)^2} \right),$$
 $$\displaystyle \int_{\O}P\delta_{(a,\lambda)}^{p}\frac{1}{\l}\frac{\partial P\delta_{(a,\lambda)}}{\partial a}=
 2\left<P\delta_{(a,\lambda)},\frac{1}{\lambda}\frac{\partial P\delta_{(a,\l)}}{\partial a}\right>+O\left(\frac{\ln(\l d)}
 {(\lambda d)^{n}}\right).$$
 \end{lem}
In the sequel, we denote
\begin{equation}\label{eijdef}\e_{ij}:=(\frac{\l_i}{\l_j}+\frac{\l_j}{\l_i}+\l_i\l_j|a_i-a_j|^2)^{(2-n)/2}.
\end{equation}
Note that, the variable $\e_{ij}$ comes from the scalar product
\begin{equation}\label{e_ij}
\displaystyle\int_{\mathbb{R}^n}\nabla\delta_{(a_i,\lambda_i)}.
\nabla\delta_{(a_j,\lambda_j)}
=\int_{\mathbb{R}^n}\delta_{(a_i,\lambda_i)}^{p}
\delta_{(a_j,\lambda_j)} =O(\e_{ij})\hbox{ for }i \neq j
,\end{equation} (see \cite{B1} page 4).
This condition means that the "interaction effect" between the
$P\d_{(a_i,\lambda_i)}$'s is negligible.\\For simplicity we shall
write $\d_i$ for $\delta_{(a_i,\lambda_i)}$ and $P\d_i$ for
$P\delta_{(a_i,\lambda_i)}$.
\begin{lem}\label{R_1}
Let $a_i,a_j\in\Omega$
 and $\lambda_i,\lambda_j > 0$ be such that $\lambda_k d_k := \lambda_k d(a_k,
 \partial\Omega)$ is very large for $k=i,j$ and $\e_{ij}$ is very small. There hold
 $$\left<P\delta_i, P\delta_{j}\right>=\overline{c}_1\left(\e_{ij}-
 \frac{H(a_i,a_j)}{(\l_i\lambda_j)^{(n-2)/2}}\right)+O(R_1),$$
 $$\left<P\delta_{j},\lambda_i\frac{\partial P\delta_{i}}{\partial \l_i}\right>=
 \overline{c}_1\left(\l_i\frac{\partial\e_{ij}}{\l_i}+\frac{n-2}{2}
 \frac{H(a_i,a_j)}{(\l_i\lambda_j)^{(n-2)/2}}\right)+O(R_1),$$
 $$\left<P\delta_{j},\frac{1}{\lambda_i}\frac{\partial P\delta_{i}}{\partial a_i}\right>=
\frac{\overline{c}_1}{\l_i}\left(\frac{\partial\e_{ij}}{\partial
a_i}-
 \frac{1}{(\l_i\lambda_j)^{(n-2)/2}}\frac{\partial H}{\partial a}(a_i,a_j)\right)+O(R_2),$$
where $\frac{\partial}{\partial a}$ (resp. $\frac{\partial}{\partial
b}$) denotes the derivative with respect to the first (resp. second)
variable of the function $(a,b)\rightarrow H(a,b)$,
$$R_1= \displaystyle \sum_{k=i,j} \frac{\ln (\l_k d_k)}{(\l_k d_k)^n}+\e_{ij}^\frac{n}{n-2}
\ln (\e_{ij}^{-1})\quad\hbox{and}\quad R_2=\sum_{k=i,j} \frac{\ln
(\l_k d_k)}{(\l_k d_k)^n}+\l_j|a_i-a_j|\e_{ij}^\frac{n+1}{n-2}.$$
\end{lem}
 \begin{lem}\label{estf}
Let $a_i,a_j\in\Omega$
 and $\lambda_i,\lambda_j > 0$ be such that $\lambda_k d_k := \lambda_k d(a_k,
 \partial\Omega)$ is very large for $k=i,j$ and $\e_{ij}$ is very small. There hold
  $$\displaystyle \int_{\O}P\delta_j^{p}P\delta_i=\left<P\delta_i,P\delta_j\right>+
  O\left(R_1\right)+ (\hbox{if }n=3)\ O\left(\e^2_{ij}\ln(\e_{ij}^{-1})^\frac{2}{3}+\frac{1}{(\l_j d_j)^2}\right),$$
 $$\displaystyle \int_{\O}P\delta_j^{p}\l_i\frac{\partial P\delta_i}{\partial \l_i}=
 \left<P\delta_i,\lambda_i\frac{\partial P\delta_i}{\partial \l_i}\right>+O\left(R_1\right)+(\hbox{if }n=3)\ O\left(\e^2_{ij}\ln(\e_{ij}^{-1})^\frac{2}{3}+\frac{1}{(\l_j d_j)^2}\right),$$
\begin{eqnarray*}\displaystyle p\int_{\O}P\delta_jP\delta_i^{p-1}\l_i\frac{\partial P\delta_i}{\partial
\l_i}=
 \left<P\delta_i,\lambda_i\frac{\partial P\delta_i}{\partial \l_i}\right>+O\left(R_1\right)
 +\ (\hbox{if }n=3)\ O\left(\e^2_{ij}\ln(\e_{ij}^{-1})^\frac{2}{3}+\frac{1}{(\l_i d_i)^2}\right),
 \end{eqnarray*}
 $$\displaystyle \int_{\O}P\delta_j^{p}\frac{1}{\l_i}\frac{\partial P\delta_i}{\partial a_i}=
 \left<P\delta_j,\frac{1}{\lambda_i}\frac{\partial P\delta_i}{\partial a_i}\right>
 +O\left(R_1\right)+ (\hbox{if }n=3)\ O\left(\e^2_{ij}\ln(\e_{ij}^{-1})^\frac{2}{3}+\frac{1}{(\l_i d_i)^2}\right),$$
\begin{eqnarray*}\displaystyle
p\int_{\O}P\delta_jP\delta_i^{p-1}\frac{1}{\l_i}\frac{\partial
P\delta_i}
 {\partial a_i}=
 \left<P\delta_i,\frac{1}{\lambda_i}\frac{\partial P\delta_i}{\partial a_i}\right>+
 O\left(R_1\right)
 +\ (\hbox{if }n=3)\ O \bigg(\e^2_{ij}\ln(\e_{ij}^{-1})^\frac{2}{3}+
 \sum_{k=i,j}\frac{1}{(\l_k d_k)^2}\bigg),
 \end{eqnarray*}
 where $R_1$ is defined in Lemma \ref{R_1}.
 \end{lem}
 To use the previous lemmas, we need some estimates of the
 nonlinearity $f_\e$ and its derivatives. We start with the
 following result contained in \cite{Pistoianew}.
\begin{lem}\label{lemmaB1}
\begin{enumerate}
\item For any $\e>0$, and any $U\in \mathbb{R}$, we have $|f_\e(U)-f_0(U)|\leq \e |U|^{p}\ln\ln(e+|U|)$.
\item For $\e$ small enough, and any $U\in \mathbb{R}$,
\begin{equation}\label{B.1} |f'_\e(U)|\leq c|U|^{p-1},\end{equation} and
\begin{equation}\label{f'}|f'_\e(U)-f'_0(U)|\leq \e
|U|^{p-1}\left(p\ln\ln(e+|U|)+\frac{1}{\ln(e+|U|)}\right).
\end{equation}
\item There exists $c>0$ such that, for $\e$ small enough and any $U,V\in
\mathbb{R}$,
\begin{equation}\label{B.2} |f'_\e(U+V)-f'_\e(U)|\leq \left\{\begin{array}{ll}c(
|U|^{p-2}+|V|^{p-2})|V|\hbox{ if }n\leq 6,\\
c( |V|^{p-1}+\e|U|^{p-1})\quad \hbox{ if }n> 6.\end{array}\right.
\end{equation}
\end{enumerate}
\end{lem}
We also have 
\begin{lem}
\begin{enumerate}\item There exists $c>0$ such that, for $\e$ small enough
and any $U,V\in \mathbb{R}$
\begin{equation}\label{B'}
|f_\e(U+V)-f_\e(U)|\leq c( |U|^{p-1}+|V|^{p-1})|V|, \quad \forall
n\geq 3.
\end{equation}
\item For $\e$ small enough, and any $U\in \mathbb{R}$,
\begin{equation}\label{B''}
|f''_\e(U)|\leq c|U|^{p-2}, \quad \forall n\geq 3.
\end{equation}
\end{enumerate}
\end{lem}
\begin{pf}
\begin{enumerate}\item By the mean value theorem there exists some $t\in
(0,1)$ such that
$$f_\e(U+V)-f_\e(U)=f'_\e(U+tV)V.$$
Hence, using \eqref{B.1} and the fact that $p-1\geq 0$ we get the
desired result.
\item Recall that
$$f'_\e(U)=\frac{|U|^{p-1}}{[\ln(e+|U|)]^\e}\left(p-\frac{\e|U|}{(e+|U|)\ln(e+|U|)}\right).$$
We see that
\begin{align*}f''_\e(U)=&\frac{\e|U|^{p-2}U}{[\ln(e+|U|)]^\e}\left(\frac{|U|-e\ln(e+|U|)}{(e+|U|)^2\ln(e+|U|)^2}
\right)\\&+
\frac{|U|^{p-3}U}{[\ln(e+|U|)]^\e}\left(p-1-\frac{\e|U|}{(e+|U|)\ln(e+|U|)}\right)
\left(p-\frac{\e|U|}{(e+|U|)\ln(e+|U|)}\right). \end{align*}
\end{enumerate}
So, for $\e$ small enough, the desired result follows.
\end{pf}\\
For $\eta>0$, $m \in \mathbb{N}$ and $
(\gamma_{1},\dots,\gamma_{m})\in\{-1,1\}^m$, let us define
$$
\begin{array}{rcl} V(m,\eta)&=&\Big\{u \in H^1_0(\Omega)\ /\ \exists
a_{1},\ldots,a_{m} \ \in \O, \ \exists
\lambda_{1},\ldots,\lambda_m>\eta^{-1}, \ \exists
\alpha_1,\ldots,\alpha_m>0 \hbox{ with}
 \\
&& \parallel u-\sum_{i=1}^{m}\alpha_i
\gamma_iP\delta_{(a_i,\lambda_i)}\parallel<\eta;\ \mid \a_i-1\mid
<\eta, \ \l_id(a_i,\partial \Omega)>\eta^{-1} \ \forall i,\
\varepsilon_{ij}<\eta \ \forall i\neq j \Big\}.
\end{array}$$
We recall that we are looking for a solution of $(P_\e)$ in a small
neighbourhood of $\sum_{i=1}^{m}\alpha_i
\gamma_iP\delta_{(a_i,\lambda_i)}$. In the following, we will
investigate the gradient of the functional $I_\e$ in $V(m,\eta)$.
\begin{pro}\label{p22}Let $n\geq 3$ and $u=\sum_{j=1}^m\a_j\gamma_jP\d_j+v\in V(m,\eta)$. For each
$i\in \{1,...,m\}$, we have the following expansion
\begin{align*}\left<\n
I_\e(u),P\d_i\right>=\gamma_i\a_i\big(1-\a_i^{p-1}\big)S_n^{\frac{n}{2}}+O\Big(\e\ln(\ln\l_i)
+\sum_{j}\frac{1}{(\l_jd_j)^{n-2}}+ \sum_{j\neq
i}\e_{ij}+\|v\|\Big),
\end{align*}
where $S_n$ is the best Sobolev constant defined in \eqref{S_n},
$d_j:=d(a_j,\partial\O)$ and $\e_{ij}$ is introduced in
\eqref{eijdef}.
\end{pro}
\begin{pf} Let $\underline{u}:= \sum_{j=1}^m\a_j\gamma_jP\delta_j$. We have
\begin{align}\label{d1d}
\displaystyle\left< \n
I_\e(u),P\delta_{i}\right>=&<u,P\delta_{i}>-\int_{\O}f_\e(u)P\delta_{i}\nonumber\\
=&<u,P\delta_{i}>-\int_{\O}\left[f_\e(u)-f_\e(\underline{u})\right]P\delta_{i}-\int_{\O}
f_\e(\underline{u})P\delta_{i}\nonumber\\
=&<\underline{u},P\delta_{i}>-A-B,
\end{align}
since $v\in E_{(a,\lambda)}$ where $E_{(a,\lambda)}$ is defined in
\eqref{vorthogonal}.\\ Using \eqref{B'}, Holder's inequality and
Sobolev embedding theorem, we obtain
\begin{align}
A&=O\left(\int_{\O}(|\underline{u}|^{p-1}+|v|^{p-1})|v|P\delta_i\right)=O(\|v\|).
\end{align}
We compute $B$.
\begin{align}
B=&\int_{\O}\left[f_\e(\underline{u})-f_\e(\a_i\gamma_iP\delta_i)\right]P\delta_i+\int_{\O}f_\e(\a_i\gamma_iP\delta_i)P\delta_i\nonumber\\
=&B_1+B_2.
\end{align}
Using again \eqref{B'} and the fact that $P\d_j\leq \d_j$, we obtain
\begin{align}\label{d1e}
B_1&=O\left(\displaystyle\int_{\O}\left(\a_i^{p-1}P\delta_i^{p-1}+|\sum_{j\neq
i}\a_j\gamma_jP\delta_j|^{p-1}\right)\sum_{j\neq i}\a_jP\delta_jP\delta_i\right)\nonumber\\
&=O\left(\sum_{j\neq i}\int_{\O}\d_i^p\d_j+\d_j^p\d_i\right)\nonumber\\
 &=O\left(\sum_{j\neq
i}\e_{ij}\right),
\end{align}
where we have used the facts that $\int_{\O}\d_j^p\d_i=O(\e_{ij})$
for all $i\neq j$ and $\a_j=O(1)$ for each $j$.\\ Notice that
\begin{equation}\label{majlnln1}0\leq\ln\ln(e+\a_iP\d_i)\leq
c\ln(\ln\l_i).\end{equation} The first claim of Lemma \ref{lemmaB1}
and \eqref{majlnln1} imply
\begin{align}\label{d1f}
B_2&:=\displaystyle\int_{\O}f_\e(\a_i\gamma_iP\delta_i)P\delta_i\nonumber\\
&=\int_{\O}f_0(\a_i\gamma_iP\delta_i)P\delta_i+\int_{\O}\left[f_\e(\a_i\gamma_iP\delta_i)-f_0(\a_i\gamma_iP\delta_i)\right]P\delta_i\nonumber\\
&=\gamma_i\a_i^p\int_{\O}P\delta_i^{p+1}+O\left(\e\ln(\ln\l_i)\right)
.
\end{align}
Combining \eqref{d1d}-\eqref{d1e}, \eqref{d1f} and Lemmas
\ref{estd}-\ref{estf}, the proof of Proposition \ref{p22} follows.
\end{pf}\\
\\
Set $$\psi_{(a,\l)}^0=\l\frac{\partial \d_{(a,\l)}}{\partial
\l}\quad\hbox{ and }\quad\psi_{(a,\l)}^1=\frac{1}{\l}\frac{\partial
\d_{(a,\l)}}{\partial a}.$$
For $\ell=0,1$, we denote by $P\psi_{(a,\l)}^\ell$ the projection of
$\psi_{(a,\l)}^\ell$ onto $ H^1_0(\O)$, defined by
\begin{eqnarray*}
-\Delta P\psi_{(a,\l)}^\ell=-\D \psi_{(a,\l)}^\ell \mbox{ in } \O,
\quad P\psi_{(a,\l)}^\ell=0 \mbox{ on}\;
\partial \O.
\end{eqnarray*}
Thus, we have
$$P\psi_{(a,\l)}^0=\l\frac{\partial P\d_{(a,\l)}}{\partial
\l}\quad\hbox{ and }\quad
P\psi_{(a,\l)}^1=\frac{1}{\l}\frac{\partial P\d_{(a,\l)}}{\partial
a}.$$
For simplicity, we denote
$\psi_{(a_i,\l_i)}^\ell$ by $\psi_{i}^\ell$ for $\ell =0,1$ and for
each $i$. Recall that simple computations show that
\begin{equation}\label{majdd}
|P\psi_i^\ell|\leq cP\d_i\leq c \d_i \hbox{ and } |\psi_i^\ell| \leq
c \d_i\hbox{ for }\ell=0,1.
 \end{equation}
To give the asymptotic expansions of the scalar product of $\nabla
I_\e(u)$ with $\l_i\frac{\partial P\d_i}{\partial \l_i}$ and
$\frac{1}{\l_i}\frac{\partial P\d_i}{\partial a_i}$ respectively, we
need the following result.
\begin{lem}\label{lemmaA1}Let $\tau $ be a positive real number
small enough. We have
\begin{align}
&\displaystyle\int_{\O}\left[f_\e(\a_i\gamma_iP\delta_i)-f_0(\a_i\gamma_iP\delta_i)\right]P\psi_i^\ell\nonumber\\&=
\left\{\displaystyle \begin{array}{ll}-\gamma_i\Gamma_1\displaystyle
\frac{\e\a_i^p}{\ln
\lambda_i}+O\left(\e^2\ln(\ln\l_i)^2+\frac{\e\ln\ln\l_i}{\l_i^{(1-\tau)n/2}d_i^n}+\frac{\e}{\ln(\l_i)^2}
+\frac{\e\ln\ln\l_i}{(\l_id_i)^{n-2}}\right),
 \hbox{ if }\ell =0; \\
\displaystyle
O\left(\e^2\ln(\ln\l_i)^2+\frac{1}{(\l_id_i)^{2n-4}}\right),
 \hbox{ if }\ell =1,\end{array}\right.\nonumber
\end{align}
where $\Gamma_1=\displaystyle \frac{2}{n-2}\int_{\mathbb{R}^n}
\delta_{(0,1)}^p(y)\ln (\delta_{(0,1)}(y))\psi_{(0,1)}^0(y)\, dy>0$
and
$\psi_{(0,1)}^0(y)=\frac{c_0(n-2)}{2}\frac{1-|y|^2}{(1+|y|^2)^{n/2}}$.
\end{lem}
\begin{pf}
Taylor's expansion with respect to $\e$ yields,
\begin{align}\label{epsilon1}
&\displaystyle\int_{\O}\left[f_\e(\a_i\gamma_iP\delta_i)-f_0(\a_i\gamma_iP\delta_i)\right] P\psi_i^\ell\nonumber\\
&=-\gamma_i\e
\int_{\O}\a_i^pP\delta_i^p\ln\ln(e+\a_iP\delta_i)P\psi_i^\ell+O\left(
\e^2\int_{\O}P\delta_i^p[\ln\ln(e+\a_iP\delta_i)]^2P\psi_i^\ell\right)\nonumber\\
&=-\gamma_i\e
\int_{\O}\a_i^pP\delta_i^p\ln\ln(e+\a_iP\delta_i)P\psi_i^\ell+
O(\e^2\ln(\ln\l_i)^2),
\end{align}
where we have used \eqref{majlnln1}.\\
Taking account of Proposition \ref{p21} $(a)$, the mean value
theorem and \eqref{majlnln1} yield
\begin{equation}\label{devg}
P\delta_i^p\ln\ln(e+\a_iP\delta_i)=\delta_i^p\ln\ln(e+\a_i\delta_i)+O(\delta_i^{p-1}\varphi_i\ln\ln\lambda_i).
\end{equation}
Hence, by using Proposition \ref{p21}, \eqref{majdd} and
\eqref{devg}, we get
\begin{align}\label{epsilon2}
&\displaystyle\int_{\O}P\delta_i^p\ln\ln(e+\a_iP\delta_i)P\psi_i^\ell\nonumber
\\&= \int_\O\delta_i^p\ln\ln(e+\a_i\delta_i)\psi_i^\ell
+O\left(\ln\ln\lambda_i\int_\O\delta_i^{p}(\varphi_i+|\psi_i^\ell-P\psi_i^\ell|)\right)\nonumber\\
&=\int_\O\delta_i^p\ln\ln(e+\a_i\delta_i)\psi_i^\ell
+O\left(\ln\ln\lambda_i\left(\int_{B_i^c}\delta_i^{p+1}+(\|\varphi_i\|_{L^\infty(B_i)}+
\|\psi_i^\ell-P\psi_i^\ell\|_{L^\infty(B_i)})\int_{B_i}\delta_i^{p}\right)\right)
\nonumber\\
&=\int_\O\delta_i^p\ln\ln(e+\a_i\delta_i)\psi_i^\ell
+O\left(\frac{\ln\ln\l_i}{(\l_id_i)^n}+\frac{1}{\l_i^\frac{n-2}{2}d_i^{n-2}}\frac{\ln\ln\l_i}{\l_i^\frac{n-2}{2}}\right)
\nonumber\\
&=\int_\O\delta_i^p\ln\ln(e+\a_i\delta_i)\psi_i^\ell
+O\left(\frac{\ln\ln\l_i}{(\l_id_i)^{n-2}}\right),
\end{align}
where $B_i$ denotes the ball of center $a_i$ and radius
$d_i/2$.\\
For $\ell=1$, we have
 \begin{align}\label{l=1}
\displaystyle \int_\O\delta_i^p\ln\ln(e+\a_i\delta_i)\psi_i^\ell=
\int_{B_i}\delta_i^p\ln\ln(e+\a_i\delta_i)\psi_i^\ell
+\int_{B_i^c}\delta_i^p\ln\ln(e+\a_i\delta_i)\psi_i^\ell
=O\left(\frac{\ln\ln\l_i}{(\l_id_i)^{n-2}}\right).
\end{align}
since the function $\delta_i^p\ln\ln(e+\a_i\delta_i)\psi_i^\ell$ is
antisymmetric with respect to $x-a_i$ in $B_i$. Thus, the desired
result follows from
\eqref{epsilon1}, \eqref{epsilon2} and \eqref{l=1}.\\
For $\ell=0$, an improvement of Lemma B.2 of \cite{Pistoianew} is
needed for our construction.
\begin{lem}\label{lemmaB2}
For any $U>0$ and $\l$ large enough, we have
$$\displaystyle\ln\ln(e+\l^\frac{n-2}{2}U)=\ln\ln(\l^\frac{n-2}{2})+\frac{2\ln U}{(n-2)\ln \l}+
\left[\ln\left(1+\frac{2\ln(e^{1-\frac{n-2}{2}\ln
\l}+U)}{(n-2)\ln\l}\right)-\frac{2\ln U}{(n-2)\ln \l}\right]$$ and
\begin{equation}\label{convergence}
\displaystyle\lim_{\l\rightarrow\infty
}(\ln\l)^2\left[\ln\left(1+\frac{2\ln(e^{1-\frac{n-2}{2}\ln
\l}+U)}{(n-2)\ln\l}\right)-\frac{2\ln U}{(n-2)\ln
\l}\right]=-\frac{2\ln (U)^2}{(n-2)^2}.
\end{equation}
\end{lem}
\begin{pf}
We have that
\begin{align}\label{devlnln}
\displaystyle\ln\ln(e+\l^\frac{n-2}{2}U)&=\ln\left(\ln(\l^\frac{n-2}{2})+\ln(e^{1-\frac{n-2}{2}\ln\l}+U)\right)\nonumber
\\&=\ln\ln(\l^\frac{n-2}{2})+\frac{2\ln U}{(n-2)\ln \l}+
\left[\ln\left(1+\frac{2\ln(e^{1-\frac{n-2}{2}\ln
\l}+U)}{(n-2)\ln\l}\right)-\frac{2\ln U}{(n-2)\ln \l}\right].
\end{align}
Set
$$\displaystyle g(t):=\frac{1}{t^2}
\left[\ln\left(1+\frac{2t\ln(e^{1-\frac{n-2}{2t}}+U)}{(n-2)}\right)-\frac{2t\ln
U}{(n-2)}\right], \quad t>0.$$ We have
$$\displaystyle \lim_{t\rightarrow 0} g(t)=-\frac{2\ln( U)^2}{(n-2)^2}.$$
Taking $t:= \ln(\l)^{-1}$, we obtain the claim.
\end{pf}\\
Let $\Omega_{\l_i}:=\l_i(\O-a_i)$. We recall that
$$\delta_{(a_i,\l_i)}=\l_i^{\frac{n-2}{2}}\delta_{(0,1)}(\l_i(.-a_i))
\hbox { and }
\psi^0_{(a_i,\l_i)}=\l_i^{\frac{n-2}{2}}\psi^0_{(0,1)}(\l_i(.-a_i)).$$Using
Lemma \ref{lemmaB2} with taking $U=\a_i\delta_{(0,1)}$, we get
\begin{align}\label{epsilon3}
\displaystyle \int_\O\delta_i^p\ln\ln(e+\a_i\delta_i)\psi_i^0=&
\l_i^{-n+\frac{n+2}{2}+\frac{n-2}{2}}\int_{\O_{\l_i}}
\delta_{(0,1)}^p(y)\ln\ln(e+\l_i^{\frac{n-2}{2}}\a_i\delta_{(0,1)}(y))\psi_{(0,1)}^0(y)\,
dy\nonumber\\
=&\ln\ln(\l_i^\frac{n-2}{2})\int_{\O_{\l_i}}
\delta_{(0,1)}^p(y)\psi_{(0,1)}^0(y)\,
dy+\frac{2}{(n-2)\ln\l_i}\int_{\O_{\l_i}} \delta_{(0,1)}^p(y)\ln
(\a_i\delta_{(0,1)})\psi_{(0,1)}^0(y)\, dy\nonumber\\
&+\int_{\O_{\l_i}}\delta_{(0,1)}^p(y)\left[\ln\left(1+\frac{2\ln(e^{1-\frac{n-2}{2}\ln
\l_i}+\a_i\delta_{(0,1)}(y))}{(n-2)\ln\l_i}\right)-\frac{2\ln
(\a_i\delta_{(0,1)})(y)}{(n-2)\ln \l_i}\right]\psi_{(0,1)}^0(y)\nonumber\\
=&A+B+C.
\end{align}
Recall that $\int_{\mathbb{R}^n}\delta_{(0,1)}^p\psi_{(0,1)}^0= 0$
and note that $B(0,\l_id_i)\subset \O_{\l_i}$. Thus
\begin{align}
A&=\ln\ln(\l_i^\frac{n-2}{2})\int_{\mathbb{R}^n}
\delta_{(0,1)}^p(y)\psi_{(0,1)}^0(y)\,
dy-\ln\ln(\l_i^\frac{n-2}{2})\int_{\O_{\l_i}^c}
\delta_{(0,1)}^p(y)\psi_{(0,1)}^0(y)\, dy\nonumber\\
&=O\left(\ln\ln\l_i\int_{B(0,\l_id_i)^c}\d_{(0,1)}^{p+1}\right)=O\left(\frac{\ln\ln\l_i}{(\l_id_i)^n}\right).
\end{align}
Concerning $B$, we have
\begin{align}
B&=\frac{2}{(n-2)\ln\l_i}\int_{\mathbb{R}^n}
\delta_{(0,1)}^p(y)\left(\ln \a_i+\ln
\delta_{(0,1)}\right)\psi_{(0,1)}^0(y)\, dy-
\frac{2}{(n-2)\ln\l_i}\int_{\O_{\l_i}^c} \delta_{(0,1)}^p(y)\ln
(\a_i\delta_{(0,1)})\psi_{(0,1)}^0(y)\, dy\nonumber\\
&=\frac{2}{(n-2)\ln\l_i}\int_{\mathbb{R}^n} \delta_{(0,1)}^p(y)\ln
(\delta_{(0,1)}(y))\psi_{(0,1)}^0(y)\, dy+
O\left(\frac{1}{\ln\l_i}\int_{B(0,\l_id_i)^c}
\delta_{(0,1)}^p(y)(|y|^2+1)\psi_{(0,1)}^0(y)\, dy\right)\nonumber\\
&=\Gamma_1\displaystyle \frac{1}{\ln \lambda_i}+O\left(\frac{1}{\ln
\l_i(\l_id_i)^{n-2}}\right)
\end{align}
by using again $\int_{\mathbb{R}^n}\delta_{(0,1)}^p\psi_{(0,1)}^0=
0$ and
since we have $|\ln\d_{(0,1)}(y)|\leq c (|y|^2+1)$ for each $y\in \mathbb{R}^n$.\\
Lastly we prove that $C$ is a remainder term. Let $\tau$ be a
positive constant small enough. We split $C$ into two parts:
 \begin{equation}\label{C}C=\int_{\O_{\l_i^{(1-\tau)/2}}}\ldots+
 \int_{\O_{\l_i}\setminus \O_{\l_i^{(1-\tau)/2}}}\ldots=C_1+C_2.
 \end{equation}
Note that $\O_{\l_i}\setminus \O_{\l_i^{(1-\tau)/2}}\subset
B(0,\l_i{\bf{d}})\setminus
 B(0,\l_i^{(1-\tau)/2}d_i)$ where ${\bf{d}}=diam(\Omega)$ and $d_i=d(a_i,\partial \Omega)>0$. In $W:=B(0,\l_i{\bf{d}})\setminus
 B(0,\l_i^{(1-\tau)/2}d_i)$ we have
 \begin{align}
\displaystyle \left|\ln\left(1+\frac{2\ln(e^{1-\frac{n-2}{2}\ln
\l_i}+\a_i\delta_{(0,1)})}{(n-2)\ln\l_i}\right)-\frac{2\ln
(\a_i\delta_{(0,1)})}{(n-2)\ln \l_i}\right|&\leq
\ln\ln(e+\l_i^\frac{n-2}{2}\a_i\d_{(0,1)})+
\ln\ln\l_i^\frac{n-2}{2}+\left|\frac{2\ln(\a_i\d_{(0,1)})}{(n-2)\ln
\l_i}\right|\nonumber\\&\leq c\ln\ln \l_i,
 \end{align}
by using \eqref{devlnln} and since
$\frac{2\ln(\a_i\d_{(0,1)})}{(n-2)\ln \l_i}$ is bounded in $W$.
Hence,
\begin{align}\label{C2}
\displaystyle C_2=O\left(\ln\ln\l_i\int_{
 B^c(0,\l_i^{(1-\tau)/2}d_i)} \d_{(0,1)}^{p+1}
 \right)=O\left(\frac{\ln\ln\l_i}{\l_i^{(1-\tau)n/2}d_i^n}\right).
\end{align}
For each $y\in \O_{\l_i^{(1-\tau)/2}}\subset B(0,
\l_i^{(1-\tau)/2}{\bf{d}})$, we have
$$\displaystyle\frac{2\ln(e^{1-\frac{n-2}{2}\ln
\l_i}+\a_i\delta_{(0,1)}(y))}{(n-2)\ln\l_i}\geq
\frac{2\ln(e\l_i^{-\frac{n-2}{2}}+\a_i
c_0(2{\bf{d}}^2\l_i^{(1-\tau)})^{-\frac{(n-2)}{2}})}{(n-2)\ln\l_i}\rightarrow
-(1-\tau)>-1, \ \hbox {as } \l_i \rightarrow \infty.$$ We deduce,
for $\l_i$ large enough, that
\begin{equation}\label{h1}\displaystyle\frac{2\ln(e^{1-\frac{n-2}{2}\ln
\l_i}+\a_i\delta_{(0,1)}(y))}{(n-2)\ln\l_i}\geq -1+\frac{\tau}{2}>-1
\hbox{ for all }y\in B(0, \l_i^{(1-\tau)/2}{\bf{d}}).
\end{equation}
 Recall that we have
\begin{equation}\label{h2}|\ln(1+x)-x|\leq c x^2\hbox{, for each } x\in [-\sigma ,\infty)
\hbox{ where } 0<\sigma<1.\end{equation} Using \eqref{h1},
\eqref{h2} and the fact that $|\ln\d_{(0,1)}|\leq
\ln(\d_{(0,1)})^2+1$, we obtain, for each $y\in B(0,
\l_i^{(1-\tau)/2}{\bf{d}})$,
\begin{align}\label{h3}
&\left|\ln\left(1+\frac{2\ln(e^{1-\frac{n-2}{2}\ln
\l_i}+\a_i\delta_{(0,1)})}{(n-2)\ln\l_i}\right)-\frac{2\ln
(\a_i\d_{(0,1)})}{(n-2)\ln \l_i}\right|\nonumber\\
&\leq c \ \frac{\ln(e^{1-\frac{n-2}{2}\ln
\l_i}+\a_i\delta_{(0,1)})^2}{\ln(\l_i)^2}+\left|\frac{2\ln(e^{1-\frac{n-2}{2}\ln
\l_i}+\a_i\delta_{(0,1)})}{(n-2)\ln\l_i}-\frac{2\ln(\a_i\d_{(0,1)})}{(n-2)\ln\l_i}\right|\nonumber\\
&\leq c\
\frac{\left(\ln(\a_i\d_{(0,1)})+\ln\left(1+e^{1-\frac{n-2}{2}\ln
\l_i}(\a_i\d_{(0,1)})^{-1}\right)\right)^2}{\ln(\l_i)^2}+\left|\frac{2\ln\left(1+e^{1-\frac{n-2}{2}\ln
\l_i}(\a_i\d_{(0,1)})^{-1}\right)}{(n-2)\ln\l_i}\right|\nonumber\\
&\leq c\ \frac{\ln(\d_{(0,1)})^2+1}{\ln(\l_i)^2}.
\end{align}
\eqref{h3} asserts that
\begin{align}&\left|\ln(\l_i)^2\delta_{(0,1)}^p\left[\ln\left(1+\frac{2\ln(e^{1-\frac{n-2}{2}\ln
\l_i}+\a_i\delta_{(0,1)})}{(n-2)\ln\l_i}\right)-\frac{2\ln
(\a_i\delta_{(0,1)})}{(n-2)\ln
\l_i}\right]\psi_{0,1}^0{\bf{1}}_{\O_{\l_i^\frac{1-\tau}{2}}}\right|\nonumber\\&\quad\leq
c \delta_{(0,1)}^{p+1}(\ln(\d_{(0,1)})^2+1)=:h(y)\hbox{ in
}\mathbb{R}^n.\nonumber
\end{align}  Now using \eqref{convergence} and the
previous uniform bounds,
 by applying the dominated convergence theorem, and
since the function $h$ is integrable in $\mathbb{R}^n$,
 we deduce
 \begin{equation}\label{C1}C_1=\frac{1}{\ln(\l_i)^2}(c'+o(1))=O\left(\frac{1}{\ln(\l_i)^2}\right),\end{equation}
 where $|c'|\leq c\int_{\mathbb{R}^n}
 \delta_{(0,1)}^p\ln(\d_{(0,1)})^2|\psi_{(0,1)}^0|$.\\
Choosing $\tau$ small enough and combining \eqref{epsilon1},
\eqref{epsilon2}, \eqref{epsilon3}-\eqref{C}, \eqref{C2} and
\eqref{C1}, the desired result follows for $\ell=1$.

\end{pf}
\begin{pro}\label{p23}
Let $n\geq 3$ and $u=\sum_{j=1}^m\a_j\gamma_j P\d_j+v\in V(m,\eta)$.
For each $i\in \{1,...,m\}$, we have the following expansion
\begin{align*}
 &\left< \nabla I_\e(u),\l_i\frac{\partial P\d_i}{\partial \l_i}\right>=\gamma_i\Gamma_1\displaystyle
\frac{\a_i^p\e}{\ln \lambda_i}+
 (n-2)\overline{c}_1\frac{\gamma_i\a_i}{2}\biggr(1-2\a_i^{p-1}\biggr)\frac{H(a_i,a_i)}{\l_i^{n-2}}
 \\&+\overline{c}_1\sum_{j\neq i}\gamma_j\a_j\Big(1 -\a_j^{p-1}-\a_i^{p-1}\Big)
 \biggl(\l_i\frac{\partial \e_{ij}}{\partial
\l_i}+\frac{n-2}{2}\frac{H(a_i,a_j)}{(\l_i\l_j)^{(n-2)/2}}
\biggr)\\&
+O\left(\frac{\e}{\ln(\l_i)^2}+\frac{1}{(\l_id_i)^{2n-4}}+\frac{1}{\l_i^{(1-\tau)n}d_i^{2n}}+\sum_{j=1}^m\e^2\ln(\ln\l_j)^2
+\sum_{j=1}^m \frac{\ln(\l_j d_j)}{(\l_jd_j)^{n}}+(\hbox{if
}n=3)\sum_{j=1}^m \frac{1}{(\l_jd_j)^{2}}\right)\\&+O\left(
\sum_{j\neq i}\e_{ij}^{\frac{n}{n-2}}\ln(\e_{ij}^{-1})+\sum_{j\neq
i}\e_{ij}^{2}\ln(\e_{ij}^{-1})^{\frac{2(n-2)}{n}} +\|v\|^2\right),
\end{align*}
where
 $\overline{c}_1$ is defined in \eqref{c_1} and $\tau$ is a
positive constant small enough.
\end{pro}
\begin{pf}
Recall that $P\psi_i^0=\l_i\frac{\partial P\delta_i}{\partial \l_i}$
and $\underline{u}:= \sum_{j=1}^m\a_j\gamma_jP\delta_j$. We have
\begin{align}\label{d2d}
\displaystyle\left< \n
I_\e(u),P\psi_{i}^0\right>=&<u,P\psi_i^0>-\int_{\O}f_\e(u)P\psi_i^0\nonumber\\
=&<u,P\psi_i^0>-\int_{\O}\left[f_\e(u)-f_\e(\underline{u})-f'_\e(\underline{u})v\right]P\psi_i^0-\int_{\O}
f_\e(\underline{u})P\psi_i^0-\int_{\O}f'_\e(\underline{u})vP\psi_i^0\nonumber\\
=&<\underline{u},P\psi_i^0>-A-B-C
\end{align}
since $v\in E_{(a,\lambda)}$.\\ By the mean value theorem, there
exists $\theta=\theta(x)\in(0, 1)$ such that
$$A=\int_{\O}\left[f_\e'(\underline{u}+\theta v)-f'_\e(\underline{u})\right]vP\psi_i^0.$$
Using \eqref{B.2}, \eqref{majdd}, Holder's inequality and Sobolev
embedding theorem, we obtain
\begin{align}
|A|&\leq c\left\{\begin{array}{ll} \int_{\O}(|\underline{u}|^{p-2}+| v|^{p-2})v^2\d_i\hbox{ if }n\leq 6,\\
 \int_{\O}(| v|^{p-1}+\e|\underline{u}|^{p-1})vP\d_i \hbox{ if }n>
 6,\end{array}\right.\nonumber\\
 &\leq c\left\{\begin{array}{ll}\|v\|^2\quad \quad\ \quad\hbox{ if }n\leq 6,\\
 \| v\|^{p}+\e \|v\| \hbox{ if }n>
 6.\end{array}\right.\nonumber
 \end{align}
 We mention that the remainder term $\|v\|^p$ is bad when $n>6$. Therefore, we need to be more precise for these dimensions. In fact, we devide the set $\O$ into three subsets and in each one, we will ameliorate this estimate. \\
 Let $\O_1 := \{ x: \a_i P\d_i \geq 2 | \sum_{j \neq i} \a_j \g_j P\d_j | \}$. Note that,
 in $\O\setminus \O_1$, it holds that $P\d_i \leq c \sum_{j\neq i} P\d_j$. Hence we get
 $$ |A| \leq c  \int_{\O\setminus \O_1}( | v|^{p-1}+\e |\underline{u}|^{p-1})
 |v | P\d_i \leq c \sum_{j\neq i} \int ( | v |^p + \e  |\underline{u}|^{p-1} |v | )
 \sqrt{ \d_i \d_j}  \leq( \| v \|^p + \e \|v\| ) \e_{ij}^{1/2} \ln(\e_{ij}^{-1})^{(n-2)/(2n)} .$$
For the integral over $\O_1$, we point out that in this set, it
holds: $\a_i P\d_i/2 \leq | \underline{u} | \leq 3 \a_i P\d_i /2$.
Let $\O_v:=\{ x:   | \underline{u} | \leq  2 | v | \}$, in $\O_1
\cap \O_v$, it holds: $P\d_i \leq c | v |$ and therefore
$$ | A | \leq c  \int_{\O_1 \cap \O_v}(| v|^{p-1}+\e|\underline{u}|^{p-1})|v| P\d_i \leq c  \int | v|^{p+1}  \leq c \| v  \| ^{p+1} .$$
Finally, for the last one, which is integrating over $\O_1 \cap (
\O\setminus \O_v) =: \O_2$,  observe that, for $x\in \O_2$, it
holds: $ |\underline{u}| \leq 2 |\underline{u} + \Theta v |$ for
each $\Theta \in (0,1)$. Thus, using the mean value theorem,
\eqref{B''}, \eqref{majdd}, Holder's inequality and the fact that
$p-2<0$ for $n>6$, we obtain
\begin{align*}
 | A | & = \Big| \int_{  \O_2 } \left[f_\e'(\underline{u}+\theta v)-f'_\e(\underline{u})\right]vP\psi_i^0 \Big|
 = \Big| \int_{  \O_2 }  f''_\e ( \underline{u} + \theta _1 \theta v ) \theta v^2 P\psi_i^0 \Big|, \hbox{ for some }\theta_1\in(0,1) \\
 &  \leq c \int _{\O_2} |  \underline{u} + \theta _1 \theta v  |^{p-2} v^2 P\d_i \leq c \int _{\O_2} |  \underline{u}   |^{p-2} v^2 P\d_i \leq c \int _{\O_2}  v^2 P\d_i ^{p-1} \leq c \| v \|^2. \end{align*}
  This completes the estimate of $A$ for $n \geq 7$ and we get
  $$ | A | \leq c \| v \| ^2 +   ( \| v \|^p + \e \|v\| ) \e_{ij}^{1/2} \ln(\e_{ij}^{-1})^{(n-2)/(2n)},
  \quad \mbox{ for } n \geq 7.$$
  Lastly, we deduce
  \begin{equation}
  |A|\leq c\left\{\begin{array}{ll}\|v\|^2\quad \quad\ \quad \quad \quad\hbox{ if }n\leq 6,\\
 \| v\|^{2}+\e^2+ \e_{ij}^{n/(n-2)} \ln(\e_{ij}^{-1})\hbox{ if }n>6
 .\end{array}\right.
  \end{equation}
To compute $B$, we write
\begin{align}
B=&\int_{\O}\left[f_\e(\underline{u})-f_\e(\a_i\gamma_iP\delta_i)-f_\e(\sum_{j\neq
i}\a_j\gamma_jP\delta_j)-f'_\e(\a_i\gamma_iP\delta_i)\sum_{j\neq
i}\a_j\gamma_jP\delta_j\right]P\psi_i^0+\int_{\O}f_\e(\a_i\gamma_iP\delta_i)P\psi_i^0\nonumber\\
&+\int_{\O}f_\e(\sum_{j\neq
i}\a_j\gamma_jP\delta_j)P\psi_i^0+\int_{\O}f'_\e(\a_i\gamma_iP\delta_i)(\sum_{j\neq
i}\a_j\gamma_jP\delta_j)P\psi_i^0\nonumber\\
=:&B_1+B_2+B_3+B_4.
\end{align}
We start with the last integral. Using \eqref{f'}, \eqref{e_ij} and
\eqref{majlnln1}, we get
\begin{align}
B_4&:=\displaystyle\int_{\O}f'_\e(\a_i\gamma_iP\delta_i)(\sum_{j\neq
i}\a_j\gamma_jP\delta_j)P\psi_i^0\nonumber\\
&=\int_{\O}f'_0(\a_i\gamma_iP\delta_i)(\sum_{j\neq
i}\a_j\gamma_jP\delta_j)P\psi_i^0+\int_{\O}\left[f'_\e(\a_i\gamma_iP\delta_i)-f'_0(\a_i\gamma_iP\delta_i)\right]
(\sum_{j\neq
i}\a_j\gamma_jP\delta_j)P\psi_i^0\nonumber\\
&=p\int_{\O}(\a_iP\delta_i)^{p-1}(\sum_{j\neq
i}\a_j\gamma_jP\delta_j)P\psi_i^0+O\left(\e\ln(\ln \l_i)\sum_{j\neq
i}\int_{\Omega}\delta_i^p\delta_j\right)\nonumber\\
&=\sum_{j\neq
i}\gamma_j\a_i^{p-1}\a_jp\int_{\O}P\delta_i^{p-1}P\delta_jP\psi_i^0+O\left(\e\ln(\ln
\l_i)\sum_{j\neq i}\e_{ij}\right).
\end{align}
Using the first claim of Lemma \ref{lemmaB1}, \eqref{majdd} and the
fact that $\ln\ln(e+|\sum_{j\neq i}\a_j\gamma_jP\d_j |)=O(
\ln\ln\l_{\max})$ where $\l_{\max}:=\max(\l_1,\ldots,\l_m)$, we have
\begin{align}
B_3&:=\displaystyle\int_{\O}f_\e(\sum_{j\neq
i}\a_j\gamma_jP\delta_j)P\psi_i^0\nonumber\\
&=\int_{\O}f_0(\sum_{j\neq
i}\a_j\gamma_jP\delta_j)P\psi_i^0+\int_{\O}\left[f_\e(\sum_{j\neq
i}\a_j\gamma_jP\delta_j)-f_0(\sum_{j\neq i}\a_j\gamma_jP\delta_j)\right]P\psi_i^0\nonumber\\
&=\int_{\O}f_0(\sum_{j\neq
i}\a_j\gamma_jP\delta_j)P\psi_i^0+O\left(\e\ln(\ln\l_{\max})\sum_{j\neq i}\int_{\Omega}\delta_j^p\delta_i\right)\nonumber\\
&=\sum_{j\neq
i}\a_j^p\gamma_j\int_{\O}P\delta_j^pP\psi_i^0+O\left(\sum_{k\neq
j;k,j\neq
i}\int_{\O}P\delta_j^{p-1}\inf(P\delta_j,P\delta_k)\delta_i+\e\ln(\ln
\l_{\max})\sum_{j\neq i}\e_{ij}\right)\nonumber\\
&=\sum_{j\neq
i}\a_j^p\gamma_j\int_{\O}P\delta_j^pP\psi_i^0+O\left(\sum_{k\neq
j}\e_{kj}^\frac{n}{n-2}\ln \e_{kj}^{-1} +\e^2\ln(\ln
\l_{\max})^2+\sum_{j\neq i}\e_{ij}^2\right).
\end{align}
Now, Lemma \ref{lemmaA1} implies that
\begin{align}
B_2&:=\displaystyle\int_{\O}f_\e(\a_i\gamma_iP\delta_i)P\psi_i^0\nonumber\\
&=\int_{\O}\left[f_\e(\a_i\gamma_iP\delta_i)-f_0(\a_i\gamma_iP\delta_i)\right]P\psi_i^0+
\int_{\O}f_0(\a_i\gamma_iP\delta_i)P\psi_i^0\nonumber\\
&=\int_{\O}\left[f_\e(\a_i\gamma_iP\delta_i)-f_0(\a_i\gamma_iP\delta_i)\right]P\psi_i^0+
\gamma_i\a_i^p\int_{\O}P\delta_i^pP\psi_i^0\nonumber\\
&=-\gamma_i\Gamma_1\displaystyle \frac{\a_i^p\e}{\ln
\lambda_i}+\gamma_i\a_i^p\int_{\O}P\d_i^pP\psi_i^0+
O\left(\e^2\ln(\ln\l_i)^2+\frac{\e}{\ln(\l_i)^2}+\frac{1}{(\l_id_i)^{2n-4}}+\frac{1}{\l_i^{(1-\tau)n}d_i^{2n}}\right).
\end{align}
In the sequel, we compute $B_1$. Let $\O_1:=\{x:\ |\sum_{j\neq
i}\a_j\gamma_jP\delta_j(x)|\leq \frac{1}{2}\a_iP\delta_i(x)\}$.
\begin{align}B_1&:=\int_{\O}\left[f_\e(\underline{u})-f_\e(\a_i\gamma_iP\delta_i)-f_\e(\sum_{j\neq
i}\a_j\gamma_jP\delta_j)-f'_\e(\a_i\gamma_iP\delta_i)\sum_{j\neq
i}\a_j\gamma_jP\delta_j\right]P\psi_i^0\nonumber \\
&=\int_{\O_1}\ldots +\int_{\O\setminus\O_1}\ldots
\nonumber \\
&=B_{11}+B_{12}.
\end{align}
Observe that, in $\O_1$, it holds $1/2 \a_iP\d_i\leq
|\a_i\gamma_iP\delta_i+\theta\sum_{j\neq
i}\a_j\gamma_jP\delta_j|\leq 3/2 \a_iP\d_i$ for each $\theta\in
(0,1)$. By using the mean value theorem and \eqref{B''}, we have
\begin{align}
|B_{11}|
&\leq\displaystyle\int_{\O_1}\left|f_\e(\underline{u})-f_\e(\a_i\gamma_iP\delta_i)
-f'_\e(\a_i\gamma_iP\delta_i)\sum_{j\neq
i}\a_j\gamma_jP\delta_j\right|\delta_i+\int_{\O_1}|f_\e(\sum_{j\neq
i}\a_j\gamma_jP\delta_j)|\delta_i\nonumber\\
&\leq
c\displaystyle\int_{\O_1}\left|f''_\e(\a_i\gamma_iP\delta_i+\theta\sum_{j\neq
i}\a_j\gamma_jP\delta_j)\right|\left(\sum_{j\neq
i}\a_jP\delta_j\right)^2\delta_i+\int_{\O_1}(\sum_{j\neq
i}\a_jP\delta_j)^p\delta_i, \hbox{ for some }\theta\in(0,1)\nonumber\\
&\leq c\displaystyle\int_{\O_1}\left(\sum_{j\neq
i}\a_jP\delta_j\right)^2\delta_i^{p-1}+c\sum_{j\neq
i}\int_{\O_1}(\delta_j \delta_i)^\frac{n}{n-2}\nonumber\\
&\leq c \sum_{j\neq i}\e_{ij}^\frac{n}{n-2}\ln (\e_{ij}^{-1})+\
(\hbox{if } n=3)\ c\sum_{j\neq i}\e_{ij}^2\ln(
\e_{ij}^{-1})^\frac{2}{3}.
\end{align}
Using again the mean value theorem, \eqref{B.1} and the fact that
$|\a_j-1|<\eta $ for each $j$, we obtain
\begin{align}
|B_{12}|
&\leq\displaystyle\int_{\O\setminus\O_1}\left(|f_\e(\a_i\gamma_iP\delta_i)|+|f'_\e(\a_i\gamma_iP\delta_i)
||\sum_{j\neq
i}\a_j\gamma_jP\delta_j|\right)P\delta_i+\int_{\O\setminus\O_1}|f_\e(\underline{u})-f_\e(\sum_{j\neq
i}\a_j\gamma_jP\delta_j)|P\delta_i\nonumber\\
&\leq
c\displaystyle\int_{\O\setminus\O_1}\a_i^pP\delta_i^{p+1}+\a_i^{p-1}P\delta_i^{p}|\sum_{j\neq
i}\a_j\gamma_jP\delta_j|+\int_{\O\setminus\O_1}|f'_\e(\sum_{j\neq
i}\a_j\gamma_jP\delta_j+\theta\a_i\gamma_iP\delta_i)|\a_iP\delta_i^2, \hbox{ for some }\theta\in(0,1)\nonumber\\
&\leq c \sum_{j\neq i}\e_{ij}^\frac{n}{n-2}\ln
(\e_{ij}^{-1})+c\int_{\O\setminus\O_1}|\sum_{j\neq
i}\a_j\gamma_jP\delta_j|^{p-1}P\delta_i^2\nonumber\\
&\leq c \sum_{j\neq i}\e_{ij}^\frac{n}{n-2}\ln (\e_{ij}^{-1})+\
(\hbox{if } n=3)\ c\sum_{j\neq i}\e_{ij}^2\ln(
\e_{ij}^{-1})^\frac{2}{3}.
\end{align}
Lastly we estimate $C$. Note that $\ln\ln(e+|\underline{u}|)=O(
\ln\ln\l_{\max})$. Using \eqref{f'} and the fact that
$$f_0'(\underline{u})=p\a_i^{p-1}P\delta_i^{p-1}+O(\sum_{j\neq
i}P\d_j^{p-1}{\bf{1}}_{\{P\d_i\leq
P\d_j\}}+P\delta_i^{p-2}P\d_j{\bf{1}}_{\{P\d_j\leq P\d_i\}}),$$ we
get
\begin{align}
C&:=\displaystyle\int_{\O}f'_\e(\underline{u})vP\psi_i^0\nonumber\\
&=\int_{\O}f'_0(\underline{u})vP\psi_i^0+\int_{\O}[f'_\e(\underline{u})-f'_0(\underline{u})]vP\psi_i^0\nonumber\\
&=p\a_i^{p-1}\int_{\O}P\delta_i^{p-1}P\psi_i^0 v+O \left(\sum_{j\neq
i}\int_{P\d_i\leq P\d_j}P\delta_j^{p-1}P\delta_i |v|
+ \int_{P\d_j\leq P\d_i}P\delta_i^{p-1}P\delta_j |v|+\|v\|\e\ln\ln\l_{\max}\right)\nonumber\\
&=C_1+O \left(\sum_{j\neq i}\int_{P\d_i\leq
P\d_j}P\delta_j^{p-1}P\delta_i |v| + \int_{P\d_j\leq
P\d_i}P\delta_i^{p-1}P\delta_j| v| +\|v\|\e\ln\ln\l_{\max}\right).
\end{align}
We need to estimate the following integral. For $k\neq j$, we have
\begin{align}
\int_{P\d_j\leq P\d_k}P\delta_k^{p-1}P\delta_j |v|
&=O\left(\|v\|\left(\int_{P\d_j\leq
P\d_k}\left(P\d_k^\frac{4}{n-2}P\d_j\right)^\frac{2n}{n+2}\right)^\frac{n+2}{2n}\right).
\nonumber
\end{align}
Observe that, for $n \geq 6$ and $k \neq j$,
$$\displaystyle\int_{P\d_j\leq P\d_k}\left(P\d_k^\frac{4}{n-2}P\d_j\right)^\frac{2n}{n+2}
\leq
\int_{\Omega}(\d_j\d_k)^\frac{n}{n-2}=O\left(\e_{kj}^\frac{n}{n-2}\ln(\e_{kj}^{-1})\right)$$
and for $n\leq 5$,
$$\displaystyle\left(\int_{P\d_j\leq P\d_k}\left(P\d_k^\frac{4}{n-2}P\d_j\right)^\frac{2n}{n+2}\right)^\frac{n+2}{2n}
\leq
\left(\int_{\O}\left(\d_k\d_j\right)^\frac{2n}{n+2}\d_k^{\frac{6-n}{n-2}\frac{2n}{n+2}}\right)^\frac{n+2}{2n}
=O\left(\e_{kj}\ln(\e_{kj}^{-1})^\frac{n-2}{n}\right).$$ Thus, we
obtain
\begin{equation}
C=C_1+O\left(\e^2\ln(\ln\l_{\max})^2+\|v\|^2+ \sum_{j\neq
i}\e_{ij}^{\frac{n+2}{n-2}}\ln(\e_{ij}^{-1})^\frac{n+2}{n}+\sum_{j\neq
i}\e_{ij}^{2}\ln(\e_{ij}^{-1})^{\frac{2(n-2)}{n}} \right).
\end{equation}
Let $B_i$ be the ball of center $a_i$ and radius $d_i/2$. Concerning
$C_1$, we have
\begin{align}\label{d2f}
C_1&=\int_{B_i}\ldots
+\int_{B_i^c}\ldots\nonumber \\
&=p\a_i^{p-1}\int_{B_i}P\delta_i^{p-1}P\psi_i^0
v+O\left(\|v\|\frac{1}{(\l_id_i)^\frac{n+2}{2}}\right)\nonumber\\
&=p\a_i^{p-1}\int_{\O}\delta_i^{p-1}\l_i\frac{\partial
\d_i}{\partial \l_i} v+O\left(\int_{B_i}\d_i^{p-1}\varphi_i
|v|+\int_{B_i}\d_i^{p-1}\l_i|\frac{\partial \varphi_i}{\partial
\l_i}|
|v|+\|v\|\frac{1}{(\l_id_i)^\frac{n+2}{2}}\right)\nonumber\\
&=O\left(\|v\|\left(\frac{1}{\l_i^\frac{n-2}{2}d_i^{n-2}}
(\int_{B_i}\d_i^\frac{8n}{n^2-4})^\frac{n+2}{2n}+\frac{1}{(\l_id_i)^\frac{n+2}{2}}\right)\right)\nonumber\\
&=O\left(\|v\|\left(\frac{1}{(\l_id_i)^\frac{n+2}{2}}+ (\hbox{if
n=6})\frac{\ln(\l_id_i)}{(\l_id_i)^4} + (\hbox{if }n<6)
\frac{1}{(\l_id_i)^{n-2}}\right)\right),
\end{align}where
we have used Proposition \ref{p21}, Holder inequality, the fact that
$v\in E_{a,\l}$ and the following computation
$$(\int_{B_i}\d_i^\frac{8n}{n^2-4})^\frac{n+2}{2n}=O\left(\frac{1}{\l_i^\frac{n-2}{2}}\ (\hbox{if }n<6)+
\frac{\ln(\l_id_i)}{\l^2}\ (\hbox{if
n=6})+\frac{d_i^{\frac{n-6}{2}}}{\l_i^2} \ (\hbox{if }n>6)\right).$$
Combining \eqref{d2d}-\eqref{d2f} and Lemmas \ref{estd}-\ref{estf},
the proof of Proposition \ref{p23} is completed.
\end{pf}
\begin{pro}\label{p24} Let $n\geq 3$ and $u=\sum_{l=1}^m\a_l\gamma_lP\d_l+v\in V(m,\eta)$. For each
 $i\in \{1,...,m\}$ and $j\in \{1,...,n\}$, we have the following expansion
\begin{align*}
&\left<\n I_\e(u),\frac{1}{\l_i}\frac{\partial P\d_i}{\partial
(a_i)_j}\right>=
 \gamma_i\biggr(\a_i^{p}-\frac{\a_i}{2}\biggr)
 \frac{\overline{c}_1}{\l_i^{n-1}}\frac{\partial H(a_i,a_i)}{
 \partial (a)_j}\\&
 +\overline{c}_1\sum_{l=1,l\neq i}^m\gamma_l\a_l\Big(1 -\a_l^{p-1}-\a_i^{p-1}\Big)
 \frac{1}{\l_i}\biggr(\frac{\partial \e_{il}}{\partial
(a_i)_j}-\frac{1}{(\l_i\l_l)^{(n-2)/2}}\frac{\partial H}{\partial
(a)_j}(a_i,a_l) \biggr)\\& +O\bigg(\sum_{l=1}^m\e^2\ln(\ln\l_l)^2
 +\sum_{l=1}^m\frac{\ln(\l_l d_l)}{(\l_l
d_l)^{n}} + (\hbox{if }n=3)\sum_{l=1}^m\frac{1}{(\l_l
d_l)^{2}}+\sum_{l\neq i}\l_l | a_i-a_l|\e_{il}^{\frac{n+1}{n-2}
}\\&+\sum_{l\neq
i}\e_{il}^{\frac{n}{n-2}}\ln(\e_{il}^{-1})+\sum_{l\neq
i}\e_{il}^{2}\ln(\e_{il}^{-1})^{\frac{2(n-2)}{n}}+\|v\|^2\bigg),
\end{align*}
where $\frac{\partial H}{\partial (a)_j}$ denotes the partial
derivative of $H$ with respect to the $j$-th component of the first
variable.
\end{pro}
The proof of Proposition \ref{p24} which we omit here is similar, up
to minor modifications, to that of Proposition \ref{p23}: arguing as
previously, taking account of Lemma \ref{lemmaA1} and using Lemmas
\ref{estd}-\ref{estf}, we obtain this expansion.
\section{Proof of Theorem \ref{t:11}}
In this section, we restrict ourselves to the case $n\geq 4$. Now
let
\begin{align*}M_{\e}=\biggr\{&(\a,\l,a,v)\in(0,+\infty)^m\times(0,+\infty)^m\times
\O_{d_0}^m\times H_0^1(\O):\big|\a_i
-1\big|<\nu_0,\\&\e\ln\ln\l_i<\nu_0,\,\forall i;
\frac{\l_i}{\l_j}<c_0, |a_i-a_j|>d_0',\forall i,j,i\neq j;\,v\in
E_{(a,\l)},\,\|v\|<\nu_0\biggr\},
\end{align*}
where $\nu_0$, $c_0$, $d_0$ and $d_0'$  are some suitable positive
constants and $\O_{d_0}=\{x\in\O/d(x,\partial\O)>d_0\}$. Let
$(\gamma_1,\ldots,\gamma_m)\in \{-1,1\}^m$, we define the function
\begin{align}  K_{\e}:M_{\e}\,\rightarrow
\R \label{g:1};\quad (\a,\l,a,v)\mapsto
I_\e\left(\sum_{i=1}^m\a_i\gamma_iP\d_{( a_i ,\l_i)}+v\right).
\end{align}
Note that the function $K_\e$ depends on the choice of
$(\gamma_1,\ldots,\gamma_m)$.\\ As described in the set $M_\e$, the
concentration speeds $\l_i$'s are comparable. For sake of
simplicity, $O(f(\l))$
 denotes any quantity dominated by $\sum^m_{i=1}f(\l_i)$. Observe
 also that the concentration points $a_i$'s are far away from the
 boundary and from each other. Hence, the parameter $\e_{ik}$ (introduced in \eqref{eijdef}) satisfies
\begin{equation}\label{eijestimate}
\e_{ik}\leq \frac{c}{(\l_i\l_k|a_i-a_k|^2)^\frac{n-2}{2}}\leq
\frac{c}{\l^{n-2}}, \ \forall i,k, \ i\neq k.
\end{equation}
\begin{pro}\label{G1}Let
$(\a,\l,a,v)\in M_{\e}$. $(\a,\l,a,v)$ is a critical point of
$K_{\e}$ if and only if $u=\sum_{i=1}^m\gamma_i\a_iP\d_i+v$ is a
critical point of $I_\e$, i.e. if and only if there exists
$\big(A,B,C\big)\in \R^m\times\R^m\times\big(\R^n\big)^m$ such that
the following holds:
\begin{align}
&\label{n1}(E_{\a_i}) \,\,\,\frac{\partial K_{\e}}{\partial
\a_i}=0,\,\ \forall
\,i\\
&\label{n2}(E_{\l_i}) \,\,\,\frac{\partial K_{\e}}{\partial
\l_i}=\mathbf{B}_i\left<\l_i\frac{\partial^2P\d_i}{\partial\l_i^2},v\right>+\sum_{j=1}^n
\mathbf{C}_{ij}\left<\frac{1}{\l_i}\frac{\partial^2P\d_{i}}{\partial
( a_i )_j\partial\l_i},v\right>,\,\,\forall\,i\\
\label{n3}&(E_{ a_i }) \,\,\,\frac{\partial K_{\e}}{\partial
 a_i }=\mathbf{B}_i\left<\l_i\frac{\partial^2P\d_{i}}{\partial\l_i\partial
 a_i },v\right>+ \sum_{j=1}^n
\mathbf{C}_{ij}\left<\frac{1}{\l_i}\frac{\partial^2P\d_{i}}{\partial
 a_i \partial ( a_i )_j},v\right>,\,\,\forall\,i\\\label{n4}
&(E_v)\,\,\,\frac{\partial K_{\e}}{\partial
v}=\sum_{i=1}^m\left(\mathbf{A}_iP\d_{i}+\mathbf{B}_i\l_i\frac{\partial
P\d_{i}}{\partial \l_i}+\sum_{j=1}^n
\mathbf{C}_{ij}\frac{1}{\l_i}\frac{\partial P\d_{i}}{\partial ( a_i
)_j}\right).
\end{align}
\end{pro}
As usual in this type of problems, we first deal with the $v$-part
of $u$. Namely, we have the following result.
\begin{pro}\label{G2}
There exists a smooth map which to any $(\e,\a,\l,a)$ verifying $\e$
is small enough and $(\a,\l,a,0)$ in $M_{\e}$, associates $\ov v \in
E_{(a,\l)},\,\|\ov v\|<\nu_0$, such that $(E_v)$ is satisfied for
some $\big(\mathbf{A},\mathbf{B},\mathbf{C}\big)\in
\R^m\times\R^m\times\big(\R^n\big)^m$. Such a $\ov v$ is unique,
minimizes $K_{\e}(\a,\l,a,v)$ with respect to $v$ in $\{v\in
E_{(a,\l)}/\|v\|<\nu_0\}$, and we have the following estimate
\begin{align}\label{vestimate}
\|\overline{v}\| \leq c\ \e\ln(\ln \l)+c \left\{
\begin{array}{ccc}
\displaystyle\frac{1}{\l^{n-2}}&\mbox{if}& n < 6,\\
\displaystyle\frac{\ln\l}{\l^{4}}&\mbox{if}& n = 6,\\
\displaystyle\frac{1}{\l^\frac{n+2}{2}}&\mbox{if}& n > 6.
\end{array}
\right.
\end{align}
\end{pro}
\begin{pf}
The proof of such a result follows the same ideas in \cite{B1} and
\cite{R}. In the following, we give a sketch of the proof. Let
$(\alpha,\lambda,a,v)\in M_\e$, $u=\sum_{i=1}^m\a_i\gamma_iP\d_i +v$
and $\underline{u}=\sum_{i=1}^m\a_i\gamma_iP\d_i $. Using \eqref{g0}
and \eqref{g:1} and expanding $K_{\e}$ with respect to $v$, we
obtain
\begin{align}\label{dev}
K_{\varepsilon}(\alpha,\lambda,a,v)=&\frac{1}{2}<u,u>-\int_{\O}F_\e(u)\nonumber\\
=&\frac{1}{2}<\underline{u},\underline{u}>
+\frac{1}{2}<v,v>-\int_{\O}\left[F_\e(u)-F_\e(\underline{u})-F'_\e(\underline{u})v-
\frac{1}{2}F_\e''(\underline{u})v^2\right]\nonumber\\
&-\int_{\O}F_\e(\underline{u})-\int_{\O}F'_\e(\underline{u})v-
\frac{1}{2}\int_{\O}F_\e''(\underline{u})v^2\nonumber\\=&K_{\varepsilon}(\alpha,\lambda,a,0)
-\int_{\O}F'_\e(\underline{u})v+\frac{1}{2}<v,v>-
\frac{1}{2}\int_{\O}F_0''(\underline{u})v^2\nonumber\\
&-
\frac{1}{2}\int_{\O}(F_\e''(\underline{u})-F_0''(\underline{u}))v^2-\int_{\O}\left[F_\e(u)
-F_\e(\underline{u})-F'_\e(\underline{u})v-
\frac{1}{2}F_\e''(\underline{u})v^2\right]\nonumber\\=&K_{\varepsilon}(\alpha,\lambda,a,0)
-\displaystyle\int_{\Omega}f_\e(\underline{u})v+\frac{1}{2}\displaystyle\int_{\Omega}\mid\nabla
v \mid^2
-\frac{n+2}{2(n-2)}\displaystyle\int_{\Omega}|\underline{u}|^{p-1}v^2+R_{\eps,\a,a,\l}(v),
\end{align}
where $R_{\eps,\a,\l,a}:=-
\frac{1}{2}\int_{\O}(F_\e''(\underline{u})-F_0''(\underline{u}))v^2-\int_{\O}\left[F_\e(u)-F_\e(\underline{u})-F'_\e(\underline{u})v-
\frac{1}{2}F_\e''(\underline{u})v^2\right]$ is a $\mathcal{C}^2$
function satisfying
$$R_{\eps,\a,\l,a}(v)=o(\|v\|^{2}), \ \  R'_{\eps,\a,\l,a}(v)=o(\|v\|), \ \ R''_{\eps,\a,\l,a}(v)=o(1)$$
uniformly with respect to $\e,\a,\l,a,\  (\a,\l,a,0)\in M_{\e}$ and $\e$ small enough.\\
Indeed, using \eqref{f'} and the fact that
$\ln\ln(e+|\underline{u}|)\leq \ln\ln\l$, it is easy to see that
\begin{align*}
\int_{\O}(F_\e''(\underline{u})-F_0''(\underline{u}))v^2=\int_{\O}(f_\e'(\underline{u})-f_0'(\underline{u}))v^2
=O\left(\e\ln\ln(\l)\int_{\O}|\underline{u}|^{p-1}v^2\right)=o(\|v\|^2)\end{align*}
 and  letting $\O_v=\{ x:      | \underline{u} | \leq 2| v |\}$, by the mean value theorem there exists $\theta =\theta(x)\in (0,1)$ such that
 \begin{align*}
\int_{\O}\left[F_\e(u)-F_\e(\underline{u})-F'_\e(\underline{u})v-
\frac{1}{2}F_\e''(\underline{u})v^2\right]
&=\frac{1}{6}\int_{\O\setminus\O_v}f_{\e}''(\underline{u}+\theta
v)v^3+O(\int_{\O_v}|v|^{p+1})\\
&=O\left(\int_{\O\setminus\O_v}|\underline{u}|^{p-2}|v|^3+\int_{\O}|v|^{p+1}\right)\\
&=O\left(\|v\|^{\inf(3,p+1)}\right)=o(\|v\|^2),\end{align*} where we
have used \eqref{B''} and the facts
that $|u|\leq 3 |v|$ in $\O_v$ and  $|F_\e(s)|\leq c |s|^{p+1}$ for each $s\in \mathbb{R}$.\\
Moreover, we know that the quadratic term in $v$, namely
$v\mapsto\int_{\Omega}\mid\nabla v \mid^2
-\frac{n+2}{n-2}\int_{\Omega}|\underline{u}|^{p-1}v^2$ is coercive,
with a modulus of coercivity bounded from below as $(\a,a,\l,0)\in
M_{\e}$ and $\e$ is sufficiently small- for a proof of this fact,
see \cite{B1,BC,R}.\\
Now we compute the linear part in $v$. Using the first point of
Lemma \ref{lemmaB1}, we get
\begin{align}\label{V1}
\displaystyle\int_{\Omega}f_\e(\underline{u})v&=
\int_{\Omega}\left[f_\e(\underline{u})-f_0(\underline{u})\right]v+\int_{\O}f_0(\underline{u})v\nonumber\\&=\int_{\O}f_0(\underline{u})v
+O\bigg(\e\int_{\O}|\underline{u}|^p\ln\ln(e+|\underline{u}|)|v|\bigg)\nonumber\\&=\int_{\O}f_0(\underline{u})v
+O\bigg(\|v\|\eps \ln\ln\l\bigg).
\end{align}
We recall that
\begin{align}\label{V2}
\displaystyle
\int_{\O}f_0(\underline{u})v=\int_{\O}|\underline{u}|^{p-1}\underline{u}v
=O\left(\bigg(\frac{1}{\l^{n-2}} \ (\hbox{if }n <
6)+\frac{\ln\l}{\l^{4}}\ (\hbox{if }n=6)
\frac{1}{\l^{\frac{n+2}{2}}} \ (\hbox{if }n > 6)\bigg)\|v\|\right).
\end{align}
The proof of the last estimate may be found in \cite{B1} and
\cite{R}. \eqref{V1} and \eqref{V2} assert that
$$\displaystyle\int_{\Omega}f_\e(\underline{u})v=O\left(\bigg(\eps
\ln\ln\l+\frac{1}{\l^{n-2}} \ (\hbox{if }n <
6)+\frac{\ln\l}{\l^{4}}\ ( \hbox{if }n=6)
+\frac{1}{\l^{\frac{n+2}{2}}}\
 (\hbox{if }n > 6)\bigg)\|v\|\right).$$
Consequently, the implicit function theorem yields the conclusion of Proposition \ref{G2}, together  with estimate \eqref{vestimate}.\\
\end{pf}\\
 As mentioned earlier, we will follow the ideas introduced in
\cite{BLR} to construct a family of solutions of $(P_{\e})$.
In fact, the result of Theorem \ref{t:11} will be obtained through a
careful analysis of \eqref{n1}-\eqref{n4} on $M_{\e}$. Once $\ov v$
is defined by Proposition \ref{G2} which we denote by $v$ for
simplicity, we estimate the corresponding numbers
$\mathbf{A},\mathbf{B},\mathbf{C}$ by taking the scalar product in
$H^1_0(\O)$ of $(E_v)$ with $P\d_i$, $\l_i{\partial P\d_i}/{\partial
\l_i}$ and $\l_i^{-1}{\partial P\d_i}/{\partial a_i }$ respectively.
Thus we get a quasi-diagonal system whose coefficients are given by
\begin{align*}&\left< P\d_i ,
P\d_j\right>=S_n^{\frac{n}{2}}\d_{ij}+O\big(\frac{1}{\l
^{n-2}}\big),\\&\left< P\d_i,\l_j \frac{\partial P\d_j}{\partial
\l_j}\right>=O\big(\frac{1}{\l^{n-2}}\big),\\& \left< P\d_i ,
\frac{1}{\l_j}\frac{\partial P\d_j}{\partial
a_j}\right>=O\big(\frac{1}{\l^{n-1}}\big),\\&\left<
\l_i\frac{\partial P\d_i}{\partial \l_i} , \l_j\frac{\partial
P\d_j}{\partial
\l_j}\right>=C_1\d_{ij}+O\big(\frac{1}{\l^{n-2}}\big),\\&
\left<\l_i\frac{\partial P\d_i}{\partial
\l_i},\frac{1}{\l_j}\frac{\partial P\d_j}{\partial
 a_j}\right>=O\big(\frac{1}{\l
^{n-1}}\big),\\&\left<\frac{1}{\l_i}\frac{\partial P\d_i}{\partial
 (a_i) _l}, \frac{1}{\l_j}\frac{\partial P\d_j}{\partial  (a_j)_h}\right>=
 C_2\d_{ij}\d_{hl}+O\big(\frac{1}{\l^{n}}\big)
\end{align*}
where $\d_{ij}$ and $\d_{hl}$ are the Kr\"{o}necker symbol and
$C_1$, $C_2$ are positive constants. \\The other hand side is given
by
\begin{align}\label{g2'}&\gamma_i\left<\frac{\partial K_{\e}}{\partial
v},P\d_i\right>=\frac{\partial K_{\e}}{\partial
\a_i};\nonumber\\&\a_i\gamma_i\left<\frac{\partial K_{\e}}{\partial
v},\l_i\frac{\partial P\d_i}{\partial
\l_i}\right>=\l_i\frac{\partial
K_{\e}}{\partial\l_i};\\
&\a_i\gamma_i\left<\frac{\partial K_{\e}}{\partial
v},\frac{1}{\l_i}\frac{\partial P\d_i}{\partial
 a_i }\right>=\frac{1}{\l_i}\frac{\partial K_{\e}}{\partial
 a_i }\nonumber.
\end{align}
Using Proposition \ref{p22} and the fact that $\frac{\partial
K_{\e}}{\partial\a_i}=\gamma_i\left<\nabla I_\e(u), P\d_i\right>$
and taking account of \eqref{vestimate},\eqref{eijestimate} and
$d_k=d(a_k,\partial \O)>d_0$ for each $k$, some computations yield
to \be\label{g3}\frac{\partial
K_{\e}}{\partial\a_i}=-(p-1)S_n^{\frac{n}{2}}\beta_i+V_{\a_i}(\e,\a,\l,a),
\ee where $$\beta:=(\beta_1,\ldots,\beta_m)=(\a_1-1,\ldots,\a_m-1)$$
and $V_{\a_i}$ is a smooth function that satisfies
\be\label{g5}V_{\a_i}(\e,\a,\l,a)=O\biggr(\beta_i^2+\e\ln(\ln
\l)+\frac{1}{\l ^{n-2}}(n\leq 5)+\frac{\ln\l}{\l ^{4}}(n=
6)+\frac{1}{\l ^{\frac{n+2}{2}}}(n\geq 7)\biggr).\ee In the same
way, from Proposition \ref{p23} we get
\be\label{g6}\l_i\frac{\partial K_{\e}}{\partial \l_i}=\Gamma_1
\frac{\e}{\ln\l_i}-\Gamma_2\left(\frac{H(a_i,a_i)}{\l_i^{n-2}}-\sum_{k\neq
i}\frac{\gamma_i\gamma_kG(a_i,a_k)}{(\l_i\l_k)^\frac{n-2}{2}}\right)+V_{\l_i}(\e,\a,\l,a),
\ee where $V_{\l_i}$ is a smooth function satisfying
\be\label{g7}V_{\l_i}=O\biggr\{\e^2\ln(\ln\l)^2+\frac{\e}{\ln(\l)^2}+\frac{1}{\l^{(1-\tau)
n}}+|\beta|(\frac{\e}{\ln \l}+\frac{1}{\l^{n-2}})\biggr\} \ee and
$$\Gamma_2= \frac{(n-2)\overline{c}_1}{2} \hbox{ where
}\overline{c}_1 \hbox{ is defined in } \eqref{c_1}.$$ Lastly, by
using Proposition \ref{p24}, we have
\be\label{g81}\frac{1}{\l_i}\frac{\partial K_{\e}}{\partial
 (a_i)_j}=\frac{\Gamma_3}{\l_i}\left(\frac{1}{\l_i^{n-2}}\frac{\partial
H}{\partial (a)_j}(a_i,a_i)-\sum_{k\neq
i}\frac{\gamma_i\gamma_k}{(\l_i\l_k)^\frac{n-2}{2}}\frac{\partial
G}{\partial (a)_j}(a_i,a_k)\right)+V_{ (a_i)_j}(\e,\a,\l,a),\ee
 where $\frac{\partial}{\partial a}$ and $\frac{\partial}{\partial b}$ denote the derivatives with respect to
 the first variable and
the second variable of the functions $(a,b)\rightarrow H(a,b)$ and
$(a,b) \rightarrow G(a,b)$, $V_{ a_i }$ is a smooth function such
that
\begin{align}\label{g9}&V_{ (a_i)_j}(\e,\a,\l,a)=
O\Big(\e^2\ln(\ln\l)^2+\frac{\ln\l}{\l^{n}}
+\frac{|\beta|}{\l^{n-1}}\Big)
\end{align}
and $$\Gamma_3=\frac{\overline{c}_1}{2}.$$ Notice that these
estimates imply
\begin{align*}&\frac{\partial K_{\e}}{\partial\a_i}=O\Big(|\beta|+\eps
\ln(\ln\l)+\frac{1}{\l^{n-2}} \ (\hbox{if }n <
6)+\frac{\ln\l}{\l^{4}}\ ( \hbox{if }n=6)
+\frac{1}{\l^{\frac{n+2}{2}}}\
 (\hbox{if }n > 6)\Big);\\
&\l_i\frac{\partial
K_{\e}}{\partial\l_i}=O\Big(\frac{\e}{\ln\l}+\frac{1}{\l^{n-2}}+\e^2\ln(\ln\l)^2\Big);\\&
\frac{1}{\l_i}\frac{\partial K_{\e}}{\partial
 (a_i)_j}=O(\frac{1}{\l^{n-1}}+\e^2\ln(\ln\l)^2).
\end{align*}
The solution of the system in $\mathbf{A}$, $\mathbf{B}$ and
$\mathbf{C}$ shows that
\begin{eqnarray}\label{g10}\begin{cases}
 \displaystyle \mathbf{A}=O\left(|\beta|+\eps
\ln(\ln\l)+\frac{1}{\l^{n-2}} \ (\hbox{if }n <
6)+\frac{\ln\l}{\l^{4}}\ ( \hbox{if }n=6)
+\frac{1}{\l^{\frac{n+2}{2}}}\
 (\hbox{if }n > 6)\right),     \\
  \displaystyle      \mathbf{B}=O\left(\frac{\e}{\ln\l}+\frac{1}{\l^{n-2}}+\e^2\ln(\ln\l)^2\right),    \\
      \displaystyle \mathbf{C}_{ij}=O\left(\frac{1}{\l^{n-1}}+\e^2\ln(\ln\l)^2+\frac{|\beta|}{\l^{n-2}}\right).
\end{cases}
\end{eqnarray}

This allows us to evaluate the right hand  side in the equations
$(E_{\l_i})$ and $(E_{ a_i })$, namely
\begin{align}
&
\mathbf{B}_i\left<\l_i\frac{\partial^2P\d_i}{\partial\l_i^2},v\right>+\sum_{j=1}^n
\mathbf{C}_{ij}\left<\frac{1}{\l_i}\frac{\partial^2P\d_i}{\partial
 (a_i)_j\partial\l_i},v\right>=O\left(\frac{1}{\l_i}\bigg(\frac{\e}{\ln\l}+\frac{1}{\l^{n-2}}+\e^2\ln(\ln\l)^2\bigg)\|v\|\right)\label{g11}\\
 & \mathbf{B}_i\left<\l_i\frac{\partial^2P\d_i}{\partial\l_i\partial  (a_i)_j},v\right>+
 \sum_{k=1}^n \mathbf{C}_{ik}\left<\frac{1}{\l_i}\frac{\partial^2P\d_i}{\partial
 (a_i)_k\partial  (a_i)_j},v\right>=O\left(\l_i\bigg(\frac{\e}{\ln\l}+
 \frac{1}{\l^{n-2}}+\e^2\ln(\ln\l)^2\bigg)\|v\|\right),\,\forall j\label{g121}
\end{align}
since
\begin{align}\label{g12'}\biggr\|\frac{\partial^2 P\d_i}{\partial
\l_i^2}\biggr\|=O\biggr(\frac{1}{\l_i^2}\biggr),\quad
\biggr\|\frac{\partial^2 P\d_i}{\partial \l_i\partial
 a_i }\biggr\|=O(1),\quad\biggr\|\frac{\partial^2 P\d_i}{\partial
 a_i ^2}\biggr\|=O(\l_i^2).
\end{align}
Now, we consider a point
$\overline{x}=(\overline{x}_1,\ldots,\overline{x}_m) \in \O^m$ such
that $\rho(\overline{x})>0$ and $\overline{x}$ is a non-degenerate
critical point of $\widetilde{{\bf{F}}}$ where
$\widetilde{{\bf{F}}}$ is introduced in \eqref{18}. We set
\begin{align*}
&\frac{1}{\l_i^{(n-2)/2}}=
\left(\frac{(n-2)\Gamma_1}{\Gamma_2}\right)^{1/2}(\Lambda_i(\overline{x})
+\zeta_i)\left(\frac{\e}{|\ln \e|}\right)^{1/2},\,\,i=1,...,m,\\
&  a_i = \overline{x}_i+\xi_i,\,\, i=1,...,m,
\end{align*}
where $\zeta_i \in \R$ and $\xi_i\in\R^N$ are assumed to be small
and
$(\L_1(x),\ldots,\L_m(x))$ is defined as the minimum of ${\bf{F}}_x$.\\
 With these changes of variables and using \eqref{g3} and
\eqref{g5}, $(E_{\a_i})$ is  equivalent to
\begin{align}
\label{g13}\beta_i=V_{\beta_i}(\e,\beta,\zeta, \xi)
=O(\e\ln|\ln\e|+|\beta|^2).
\end{align}
Now, using the changes of variables, an easy computation shows that
\begin{equation}\label{epsilonlambda}
\displaystyle\ln \l_i
=\frac{|\ln\e|}{n-2}\left(1+O\big(\frac{\ln|\ln\e|}{|\ln\e|}\big)\right)
\quad\hbox{ and }\quad\frac{1}{\ln \l_i}
=\frac{n-2}{|\ln\e|}+O\left(\frac{\ln|\ln\e|}{|\ln\e|^2}\right).
\end{equation}
Moreover, we have
\begin{equation}
H(a_i,a_i)=H(\overline{x}_i,\overline{x}_i)+2\frac{\partial
H}{\partial a}(\overline{x}_i,\overline{x}_i)\xi_i+O(|\xi_i|^2)
\end{equation}
and
\begin{equation}
G(a_i,a_j)=G(\overline{x}_i,\overline{x}_j)+\frac{\partial
G}{\partial a}(\overline{x}_i,\overline{x}_j)\xi_i+\frac{\partial
G}{\partial b}(\overline{x}_i,\overline{x}_j)\xi_j+O(|\xi|^2).
\end{equation}
Set $\overline{\L}_i=\L_i(\overline{x})$ for each $1\leq i\leq m$,
we have
\begin{equation}\label{condition}
1-H(\overline{x}_i,\overline{x}_i)\overline{\Lambda}_i^2+\displaystyle\sum_{j\neq
i}
\gamma_i\gamma_jG(\overline{x}_i,\overline{x}_j)\overline{\Lambda}_i\overline{\Lambda}_j=0
\end{equation} since $\overline{\L}=(\overline{\L}_1,\ldots,\overline{\L}_m)$
is a critical point of ${\bf{F}}_{\overline{x}}$. \\ \eqref{g6},
\eqref{g7} and \eqref{g11} imply that $(E_{\l_i})$ is equivalent,
while using \eqref{epsilonlambda}-\eqref{condition}, to
\begin{align}\label{g14}
&\left(2H(\overline{x}_i,\overline{x}_i)\overline{\Lambda}_i-\sum_{j\neq
i}\gamma_i\gamma_jG(\overline{x}_i,\overline{x}_j)\overline{\Lambda}_j\right)\zeta_i-\sum_{j\neq
i}\gamma_i\gamma_jG(\overline{x}_i,\overline{x}_j)\overline{\Lambda}_i\zeta_j\nonumber\\&
+ \left(2\frac{\partial H}{\partial
a}(\overline{x}_i,\overline{x}_i)\overline{\Lambda}_i^2-\sum_{j\neq
i}\gamma_i\gamma_j\frac{\partial G}{\partial
a}(\overline{x}_i,\overline{x}_j)\overline{\Lambda}_i\overline{\Lambda}_j\right)\xi_i-\sum_{j\neq
i}\gamma_i\gamma_j\frac{\partial G}{\partial
b}(\overline{x}_i,\overline{x}_j)\overline{\Lambda}_i\overline{\Lambda}_j\xi_j
=V_\zeta(\e,\beta,\zeta,\xi), \end{align} where
$$V_\zeta(\e,\beta,\zeta,\xi)=O\left(\frac{\ln|\ln\e|}{|\ln\e|}
+|\beta|^2+|\zeta|^2+\xi^2\right).$$
We also have
\begin{equation}\label{g811}
\frac{\partial H}{\partial a}(a_i,a_i)=\frac{\partial H}{\partial
a}(\overline{x}_i,\overline{x}_i)+\frac{\partial^2 H}{\partial
a^2}(\overline{x}_i,\overline{x}_i)\xi_i+\frac{ \partial^2
H}{\partial a\partial
b}(\overline{x}_i,\overline{x}_i)\xi_i+O(|\xi|^2)
\end{equation}
and
\begin{equation}\label{g91}
\frac{\partial G}{\partial a}(a_i,a_j)=\frac{\partial G}{\partial
a}(\overline{x}_i,\overline{x}_j)+\frac{\partial^2 G}{\partial
a^2}(\overline{x}_i,\overline{x}_j)\xi_i+\frac{ \partial^2
G}{\partial a\partial
b}(\overline{x}_i,\overline{x}_j)\xi_j+O(|\xi|^2)
\end{equation}
Through \eqref{g81}, \eqref{g9} and \eqref{g121}, and by using
\eqref{g811}, \eqref{g91}, we get that $(E_{ a_i})$ is equivalent to
\begin{align}\label{g18}
&\frac{\partial H}{\partial
a}(\overline{x}_i,\overline{x}_i)\zeta_i-\sum_{j\neq
i}\gamma_i\gamma_j\frac{\partial G}{\partial
a}(\overline{x}_i,\overline{x}_j)\zeta_j+ \left(\frac{\partial^2
H}{\partial
a^2}(\overline{x}_i,\overline{x}_i)\overline{\Lambda}_i+\frac{\partial^2
H}{\partial a\partial
b}(\overline{x}_i,\overline{x}_i)\overline{\Lambda}_i-\sum_{j\neq
i}\gamma_i\gamma_j\frac{\partial^2 G}{\partial
a^2}(\overline{x}_i,\overline{x}_j)\overline{\Lambda}_j\right)\xi_i\nonumber\\&-\sum_{j\neq
i}\gamma_i\gamma_j\frac{\partial^2 G}{\partial a\partial
b}(\overline{x}_i,\overline{x}_j)\overline{\Lambda}_j\xi_j
=V_{\xi}(\e,\beta,\zeta,\xi),
\end{align}
where
$$V_{\xi}(\e,\beta,\zeta,\xi)=O\left(\frac{\ln|\ln\e|}{|\ln\e|}
+|\beta|^2+|\zeta|^2+\xi^2\right).$$ Furthermore, \eqref{g13},
\eqref{g14} and \eqref{g18} may be written as
\begin{eqnarray}\label{g16'}
\begin{cases}\beta=V(\e,\beta,\zeta,\xi),\\
L(\zeta,\xi)=W(\e,\beta,\zeta,\xi),
\end{cases}
\end{eqnarray}
where $L$ is a fixed linear operator of $\R^{m(n+1)}$ defined by
\eqref{g14} and \eqref{g18} and $V,\,W$ are smooth functions
satisfying
\begin{eqnarray}\label{g17}
\begin{cases}V(\e,\beta,\zeta,\xi)=O(\e\ln|\ln\e|+|\beta|^2),\\
W(\e,\beta,\zeta,\xi)=O\big(\frac{\ln|\ln\e|}{|\ln\e|}
+|\beta|^2+|\zeta|^2+\xi^2\big).
\end{cases}
\end{eqnarray}
Moreover, a simple computation shows that the determinant of $L$ is
proportional to the determinant of
$\widetilde{{\bf{F}}}''(\overline{x})$. $\overline{x}$ being a
non-degenerate critical point of $\widetilde{{\bf{F}}}$ by
assumption, $L$ is invertible, and Brouwer's fixed point theorem
 ensures, provided that $\e$ is small
enough, the existence of a solution $(\b_{\e},\zeta_{\e},\xi_{\e})$
to $(\ref{g16'})$, such that
$$|\b_{\e}|=O(\e\ln|\ln\e|),\quad |\zeta_{\e}|=O\left(\frac{\ln|\ln\e|}{|\ln\e|}\right)\quad \hbox{and}
\quad| \xi_{\e}|=O\left(\frac{\ln|\ln\e|}{|\ln\e|}\right).$$ By
construction, $u_{\e}=\displaystyle \sum ^{m}_{i=1}
\a_{i,\e}\gamma_i P\d_{(a_{i,\e},\l_{i,\e})}+
\overline{v}(\e,\a_{\e},\l_{\e},a_{\e}) \in H_0^1(\O)$ with
$$\a_{i,\e}=1+\b_{i,\e}; \ \ \ \ \displaystyle\frac{1}{\l_{i,\e}}=
\left(\frac{(n-2)\Gamma_1}{\Gamma_2}\right)^{\frac{1}{n-2}}(\Lambda_i(\overline{x})
+\zeta_{i,\e})^{\frac{2}{n-2}}\left(\frac{\e}{|\ln
\e|}\right)^{\frac{1}{n-2}}; \ \ \ \
a_{i,\e}=\overline{x}_i+\xi_{i,\e},$$ is a critical point of
$I_{\e}$, whence
$$-\Delta u_{\e}=\frac{|u_\e|^{p-1}u_\e}{[\ln (e+|u_\e|)]^\e} \hbox{ in
}\O.$$
 The proof of
Theorem \ref{t:11} is thereby completed.\hfill$\Box$

\end{document}